\documentclass[a4paper]{article}
\usepackage{amsmath,amsthm,amssymb}

\begin{document}

\numberwithin{equation}{section}

{\theoremstyle{definition}\newtheorem{definition}{Definition}[section]
\newtheorem{notation}[definition]{Notation}
\newtheorem{remnot}[definition]{Remarks and notation}
\newtheorem{terminology}[definition]{Terminology}
\newtheorem{remark}[definition]{Remark}
\newtheorem{remarks}[definition]{Remarks}
\newtheorem{example}[definition]{Example}
\newtheorem{examples}[definition]{Examples}}
\newtheorem{proposition}[definition]{Proposition}
\newtheorem{lemma}[definition]{Lemma}
\newtheorem{theorem}[definition]{Theorem}
\newtheorem{corollary}[definition]{Corollary}
\newcommand{\cB}{\mathcal{B}}
\newcommand{\R}{\mathbb{R}}
\newcommand{\cO}{\mathcal{O}}
\newcommand{\cF}{\mathcal{F}}
\newcommand{\cE}{\mathcal{E}}
\newcommand{\B}{\mathcal{L}}
\newcommand{\K}{\mathcal{K}}
\newcommand{\cL}{\mathcal{L}}
\newcommand{\cR}{\mathcal{R}}
\newcommand{\cW}{\mathcal{W}}
\newcommand{\cV}{\mathcal{V}}
\newcommand{\C}{\mathbb{C}}
\newcommand{\full}{_{\text{\rm f}}}
\newcommand{\red}{_{\text{\rm r}}}
\newcommand{\cA}{\mathcal{A}}
\newcommand{\cC}{\mathcal{C}}
\newcommand{\cI}{\mathcal{I}}
\newcommand{\Q}{\mathbb{Q}}
\newcommand{\Z}{\mathbb{Z}}
\newcommand{\cP}{\mathcal{P}}
\newcommand{\cU}{\mathcal{U}}
\newcommand{\cH}{\mathcal{H}}
\newcommand{\resp}{{\it resp.}\/ }
\newcommand{\id}{{\hbox{id}}}
\newcommand{\ie}{{\it i.e.}\/ }
\newcommand{\eg}{{\it e.g.}\/ }
\newcommand{\cf}{{\it cf.}\/ }
\newcommand{\Aut}{{\rm Aut}}
\newcommand{\Hom}{{\rm Hom}}
\newcommand{\gog}{\mathfrak{g}}
\newcommand{\hoh}{\mathfrak{h}}
\newcommand{\X}{\mathcal{X}}
\newcommand{\Hol}{\mathcal{H}}
\def\gpd{\,\lower1pt\hbox{$\longrightarrow$}\hskip-.24in\raise2pt
             \hbox{$\longrightarrow$}\,}

%\begin{document}

\begin{center}
{\Large\bf The holonomy groupoid of a singular foliation
\footnote{Research supported by EU RTN "Quantum spaces and Noncommutative Geometry", contract number HPRN-CT-2002-00280}
}

\bigskip

{\sc by Iakovos Androulidakis and Georges Skandalis}

\end{center}

{\footnotesize 
Institut f\"{u}r Mathematik,  Universit\"{a}t Z\"{u}rich,
\vskip-4pt Winterthurerstrasse 190,  CH-8057 Z\"{u}rich
\vskip-4pt e-mail: iakovos.androulidakis@math.unizh.ch

\vskip 2pt Institut de Math{\'e}matiques de Jussieu, UMR 7586 
\vskip -4ptCNRS - Universit\'e Diderot - Paris 7
\vskip-4pt 175, rue du Chevaleret, F--75013 Paris
\vskip-4pt e-mail: skandal@math.jussieu.fr
}
\bigskip
\everymath={\displaystyle}

\begin{abstract}\noindent
We construct the holonomy groupoid of \emph{any} singular
foliation. In the regular case this groupoid coincides with the
usual holonomy groupoid of Winkelnkemper
(\cite{Winkelnkemper:holgpd}); the same holds in the singular
cases of \cite{Pradines, Bigonnet-Pradines, Debord:2001LILA,
Debord:2001HGSF}, which from our point of view can be thought of
as being ``almost regular''.  In the general case, the holonomy
groupoid can be quite an ill behaved geometric object. On the
other hand it often has a nice longitudinal smooth structure.
Nonetheless, we use this groupoid to generalize to the singular
case Connes' construction of the $C^*$-algebra of the foliation.
We also outline the construction of a longitudinal
pseudo-differential calculus; the analytic index of a
longitudinally elliptic operator takes place in the $K$-theory of
our $C^*$-algebra.

\noindent In our construction, the key notion is that of a
\emph{bi-submersion} which plays the role of a \emph{local Lie
groupoid} defining the foliation. Our groupoid is the quotient of
germs of these bi-submersions with respect to an appropriate
equivalence relation.
\end{abstract}

\section*{Introduction}

A foliated manifold is a manifold partitioned into immersed
submanifolds (leaves). Such a partition often presents
singularities, namely the dimension of the leaves needs not be
constant. Foliations arise in an abundance of situations, e.g. the fibers of a
submersion, the orbits of a Lie group action or the structure of a
manifold with corners. They are quite important tools in geometric
mechanics; for instance, foliations appear in the problem of Hamiltonian
reduction (\cf \eg the monograph \cite{OR}). Also, every Poisson manifold is
determined by its symplectic foliation (cf \cite[pages 112-113]{CF2}).

The relationship between foliations and groupoids is very well
known: to any Lie groupoid corresponds a foliation, namely the
leaves are the orbits of the groupoid. Conversely, to a regular
foliation there corresponds its \emph{holonomy groupoid}
constructed by Ehresmann \cite{Ehresmann} and Winkelnkemper
\cite{Winkelnkemper:holgpd} (see also \cite{Bigonnet-Pradines},
\cite{Pradines:holgpd}).  The holonomy groupoid is a (not
necessarily Hausdorff) Lie groupoid whose orbits are the leaves
and which in a sense is minimal with this condition.

This holonomy groupoid has been generalized to singular cases by
various authors.   In particular, a construction suggested by
 Pradines and Bigonnet  \cite{Pradines,
Bigonnet-Pradines} was carefully analyzed, and its
precise range of applicability found, by Claire Debord \cite{Debord:2001LILA, Debord:2001HGSF}.
These authors dealt with those singular foliations that can be
defined by an {\it almost injective} Lie algebroid, namely a Lie
algebroid whose anchor map is injective in a dense open subset of
the base manifold. It was shown that every such Lie algebroid is
integrable, and the Lie groupoid that arises is the holonomy
groupoid of the foliation, in the sense that it is the smallest
among all  Lie groupoids which realize the
foliation. Note that the integrability of almost injective Lie
algebroids was reproved as a consequence of Crainic and Fernandes'
characterization of integrable Lie algebroids in
\cite{Crainic-Fernandes}.

Both in the regular case and the singular case studied by
Bigonnet, Pradines and Debord, the holonomy groupoid was required
to be:
\begin{enumerate}
\item {\it Lie}, namely the space of arrows to have a smooth
structure such that the source and target maps are surjective
submersions and all the other structure maps (multiplication,
inversion, unit inclusion) to be appropriately smooth. \item {\it
minimal}, namely every other Lie groupoid which
realizes the foliation  has an open subgroupoid which maps onto the holonomy groupoid by
a smooth morphism of Lie groupoids.
\end{enumerate}
Minimality is essential mainly because it does away with
unnecessary isotropy and this ensures that the holonomy groupoid
truly records all the information of the foliation. On the other
hand, the smoothness requirement leaves room for discussion. As we will, see, in
order to construct the $C^*$-algebra it is enough that
the $s$-fibers be smooth, and we can even get around with this requirement. In fact, smoothness of the full
space of arrows is exactly what prevents generalizing their
construction to a singular foliation which comes from an arbitrary
Lie algebroid.

This paper addresses the problem of the existence of a holonomy
groupoid for a singular foliation without assuming the foliation
to be defined by some algebroid. The properties of this groupoid may be summarized in
the following theorem:
\begin{theorem}
Let $\cF$ be a (possibly singular) Stefan-Sussmann foliation on a
manifold $M$. Then there exists a topological groupoid $\Hol(\cF)
\gpd M$ such that:
\begin{itemize}
\item Its orbits are the leaves of the given foliation $\cF$.
\item $\Hol(\cF)$ is minimal in the sense that if $G \gpd M$ is a Lie groupoid which defines the foliation $\cF$ then
there exists an open subgroupoid $G_0$ of $G$ (namely its $s$-connected component) and a morphism of groupoids $G_0 \to \Hol(\cF)$ which is onto. 
\item If $\mathcal{F}$ is
regular or almost regular \ie defined from an almost injective Lie
algebroid, then $\Hol(\cF) \gpd M$ is the holonomy groupoid given
in \cite{Winkelnkemper:holgpd, Pradines:holgpd, Pradines,
Bigonnet-Pradines, Debord:2001LILA, Debord:2001HGSF}.
\end{itemize}
\end{theorem}

Note that the only possible topology of  $\Hol(\cF)$ is quite pathological in the non quasi-regular case.  On the other hand its source-fibers are often smooth manifolds.

But before giving an overview of our approach, let us first
clarify the way a singular foliation is understood. The
established definition of a foliation, following the work of
Stefan \cite{Stefan} and Sussmann \cite{Sussmann} is the partition
of a manifold $M$ to the integral submanifolds of a (locally)
finitely generated module of vector fields  which is integrable \ie stable under
the Lie bracket. In the regular case there is only one possible
choice for this module: this is the module of vectors tangent
along the leaves, it forms a vector subbundle and a Lie
subalgebroid of $TM$. In the singular case though, the partition
of $M$ into leaves no longer defines uniquely the module of vector
fields. Take for example the foliation of $\R$ into three leaves:
$\R^{-}, \{ 0 \}$ and $\R^{+}$. These can be regarded as integral
manifolds of any of the vector fields $X_{n} = x^{n}
\frac{\partial}{\partial x}$, yet the modules generated by each
$X_{n}$ are genuinely different.

Our construction depends on the choice of the prescribed submodule
of vector fields defining the singular foliation. Actually, we call \emph{foliation}  precisely such an integrable (locally) finitely generated submodule of vector fields.

Note that the almost regular case is the case where the prescribed module of vector fields is isomorphic to the sections of a vector bundle - \ie is locally free.

The construction of the holonomy groupoid by Winkelnkemper in the
regular case is based on the notion of the holonomy of a path.
The point of view of Bigonnet-Pradines is a little different: they
first look at abstract holonomies as being local diffeomorphisms
of local transversals, and obtain in this way a ``big'' groupoid.
The holonomy groupoid is just the $s$-connected component of this
big groupoid - or equivalently the smallest open subgroupoid \ie the subgroupoid generated by a
suitable open neighborhood of the space of units.

In a sense, our construction follows the ideas of
Bigonnet-Pradines. So we will construct a big groupoid, as well as
`a suitable open neighborhood of the space of units'.
\begin{itemize}
\item Even without knowing what the (big) holonomy groupoid might
be, it is quite easy to define what a submersion to this groupoid
could mean: this is our notion of  {\em bi-submersions}. A
bi-submersion is understood merely as a manifold together with two
submersions over $M$ which generate (locally) the foliation
$\mathcal{F}$ in the same way as a Lie groupoid would. The
composition (pulling back) and inversion of bi-submersions are
quite easily defined. Let us emphasize the fact that two different bi-submersions may very well be of different dimension. On the other hand, there is a natural criterion for two
points of various bi-submersions to have the same image in the
holonomy groupoid: this is the notion of equivalence of two (germs)
bi-submersions at a point. The (big) holonomy groupoid is thus
defined as being the quotient of the union of all possible
bi-submersions by this equivalence relation.
 \item Exponentiating
small vector fields defining the foliation gives rise to such a
bi-submersion. The image of this bi-submersion is the desired
`suitable open neighborhood of the space of units' in the big
groupoid. The (path) holonomy groupoid is thus the subgroupoid
generated by this neighborhood.
\end{itemize}

The ``big'' holonomy groupoid takes into account all
possible holonomies and the path holonomy groupoid is the smallest
possible holonomy groupoid. In order to treat them simultaneously as well as all possible intermediate groupoids, we introduce a notion of \emph{atlas} of
bi-submersions: loosely speaking, an atlas is a family of
bi-submersions which is stable up to equivalence by composition
and inverse. The quotient of this atlas by the equivalence
relation given above is algebraically a groupoid over $M$. The
notion of atlas of bi-submersions we introduce here seems to be in
the spirit of Pradine's ideas in \cite{Pradines:2004}.

Unless the foliation is almost regular,  the (path) holonomy groupoid (endowed with the
quotient topology) is topologically {\tt quite an ugly space. Its topology is by no means Hausdorff: a sequence of points may converge to \emph{every point of a submanifold}.} On the other hand, we show that
when restricted over a leaf $L$ with the leaf topology, it often carries
a natural smooth structure making the source (and target) map a
smooth submersion. When this happens, this restriction is actually a principal
groupoid. The stabilizer of a point \ie the holonomy group is a
Lie group, whose Lie algebra is very naturally expressed in terms
of the foliation $\mathcal{F}$. Every source (target) fiber has
the same dimension as the underlying holonomy cover of the leaf.
Of course, this dimension is not the
same for different leaves - unless the foliation is almost regular.

\medskip In the (quasi) regular case, the holonomy groupoid  of a foliation is the first
step towards the construction of the $C^*$-algebra of the foliation which allows to obtain a number of different important results:
\begin{itemize}
\item Important information of any foliation lies in the space of
leaves, namely the quotient of the manifold by the equivalence
relation of belonging in the same leaf. This space presents
considerable topological pathology, even in the regular case. A.
Connes showed in \cite{Connes}  that it is appropriately replaced by the foliation $C^*$-algebra constructed from the holonomy groupoid.
\item Using the holonomy groupoid and its $C^*$-algebra, A. Connes
(and partly the  second author) \cite{Connes0, Connes1, Connes, Connes-Skandalis} developed a longitudinal pseudodifferential calculus on foliations and an index
theory for foliations. \item In another direction, by extending
the construction of the $C^*$-algebra to an arbitrary Lie groupoid
we get a formal deformation quantization of the Poisson structure
on the dual of an integrable Lie algebroid (see \cite{Landsman}
and \cite{Landsman-Ramazan}). \item  Furthermore, the operators
corresponding to the $C^*$-algebra of a Lie groupoid are families
of pseudodifferential operators on the source fibers of the
groupoid in hand, defined by Nistor, Weinstein and Xu
\cite{Nistor-Weinstein-Xu} and independently by Monthubert and
Pierrot \cite{Monthubert-Pierrot}. In the case of the holonomy
groupoid such operators play an important role in both index
theory and deformation quantization.
\end{itemize}

Connes' construction of the $C^*$-algebra of the foliation  involves working with the algebra of smooth
functions on the arrow space of the holonomy groupoid (first
defining involution and convolution and then completing it
appropriately). Therefore it can no longer be applied for this
holonomy groupoid, since the topology of its arrow space is very
pathological, as a quotient topology, and the functions on it are
highly non-continuous.

Nevertheless, we may still define some $C^*$-algebra(s). In fact
we just push one small step further Connes' idea for the
construction of the $C^*$-algebra of the foliation in the case the
groupoid is not Hausdorff (\cf \cite{Connes1, Connes}; see also
\cite{KhoshSkand}).  Every bi-submersion is a smooth manifold, so
we may work with the functions defined on it instead. It is then
quite easy to decide when two functions defined on different
bi-submersions should be identified, \ie we define a suitable
quotient of the vector space spanned by smooth compactly supported
functions on bi-submersions. The inversion and composition
operations of the atlas are translated to an involution and
convolution  at the quotient level, making it a $*$-algebra. One
has to be somewhat careful with the formulation of this, in order
to keep track of the correct measures involved in the various
integrations.  For this reason we work with suitable
half-densities (following Connes \cf \cite{Connes}). It is finally
quite easy to produce an $L^1$-estimate which allows us to
construct the full $C^*$-algebra.
Actually, we describe all the representations of this $C^*$-algebra as being the integrated form of representations of the groupoid, pretty much in the spirit of \cite{Renault} (\cf also \cite{FackSkand}).

To construct the reduced $C^*$-algebra, we have to choose a
suitable family of representations which should be the `regular'
ones. In `good' cases, our groupoid carries a natural longitudinal
smooth structure which allows us to define the associated family
of regular representations. An alternate set of representations
which exist in all cases and which can be thought in a sense as
the regular ones are the natural representations on the space of
$L^2$-functions on a leaf. To treat all possibilities in the same
frame, we introduce a notion of a holonomy pair which consists of
an atlas together with a longitudinally smooth quotient of the associated holonomy groupoid. The
associated regular representations are those on $L^2(G_x)$ where
$x\in M$.

With these in hand we outline the analytic index map following by
and large the pattern of the tangent groupoid introduced by A.
Connes (\cf \cite{Connes}).

\bigskip
The text is structured as follows: \begin{itemize}
\item Section 1 is an account of well known fundamental facts about singular foliations and the various groups and  pseudogroups of (local) diffeomorphisms which respect them. In particular, we recall how to a Stefan-Sussmann foliation is associated the partition into leaves and discuss the longitudinal smooth structure and topology. 
\item  In Section 2 we study
bi-submersions and their bisections. We discuss
the composition of bi-submersions, equivalence of germs of bi-submersions and those bi-submersions which have bisections inducing
the identity map.
\item  In Section 3 we construct the holonomy groupoid of a foliation. We first define  the notion of an atlas of  bi-submersions, discuss several possible atlases, and then construct the groupoid of an atlas.
\item In section 4 we construct the reduced $C^*$-algebras of the foliation, 
\item In section 5 we define the representations of the holonomy groupoid and show that they are in a   one to one correspondence with the representations of (the full) $C^*$-algebra of the foliation.
\item We conclude with a
short discussion about further developments and in particular the
analytic index map (Section 6).
\end{itemize}

{\tt Actually, this pseudodifferential calculus is in a sense our main motivation for generalizing the $C^*$-algebra of a foliation to the singular case. Indeed, as in the regular case, our $C^*$-algebra is designed to deal with longitudinal differential operators that are elliptic along the foliation:
\begin{itemize}
\item  it is a receptacle for resolvants of such operators
\item its $K$-theory is a receptacle for index problems.
\end{itemize}
}

In a forthcoming paper we intend to give the complete formulation of the
pseudodifferential calculus and its correspondence with the
$C^*$-algebras given here. We also intend to extend our
construction to a case of `longitudinally smooth foliations' in
the spirit of continuous family groupoids of Paterson.

\section*{Acknowledgements}

The first author would like to thank Claire Debord and Alberto
Cattaneo for several illuminating discussions. Also the Operator
Algebras Research group at Institute Math{\'e}matiques de Jussieu
and the Research group under Alberto Cattaneo in the University of
Zurich for various comments and insight along the course of this
research.

\section{Preliminaries} \label{sec.prelim}

In this section, we recall the definitions and the results on
foliations that will be used in the sequel. For the reader's
convenience we include (or just sketch) most proofs.

\subsection{Bundles and submodules}

Let $M$ be a manifold and let $E$ be a smooth real vector bundle
over $M$.

\begin{enumerate}
\item Denote by $C^\infty(M;\R)$ or just $C^\infty(M)$ the algebra
of smooth real valued functions on $M$. Denote by
$C_c^{\infty}(M)$ the ideal in $C^\infty(M)$ consisting of smooth
functions with compact support in $M$.

\item Denote by $C^\infty(M;E)$  the $C^\infty(M)$-module of
smooth sections of the bundle $E$. Denote by $C_c^{\infty}(M;E)$
the submodule of $C^\infty(M;E)$ consisting of smooth sections
with compact support in $M$.

\item Let  $\cE$ be a submodule of $C_c^\infty(M;E)$. The
submodule $\widehat \cE\subset C^\infty(M,E)$ of \emph{global
sections} of $\cE$ is the set of sections $\xi \in C^\infty(M;E)$
such that, for all $f\in C_c^\infty(M)$, we have $f\xi \in \cE$.

The module $\cE$ is said to be \emph{finitely generated} if there
exist global sections $\xi_1,\ldots,\xi_n$ of $\cE$ such that
$\cE=C_c^\infty(M)\xi_1+\ldots +C_c^\infty(M)\xi_n$.

\item Let $N$ be a manifold and $p:N\to M$ be a smooth map. Denote
by $p^*(E)$ the pull-back bundle on $N$. If $\cE$ is a submodule
of $C_c^\infty(M;E)$, the \emph{pull-back module} $p^*(\cE)$ is
the submodule of $C_c^\infty (N;p^*(E))$ generated by $f\,
(\xi\circ p)$, with $f\in C_c^\infty (N)$ and $\xi\in \cE$.

If $N$ is a submanifold of $M$, the module $p^*(\cE)$ is called
the \emph{restriction} of $\cE$ to $N$.

\item A submodule $\cE$ of $C_c^\infty(M;E)$ is said to be
\emph{locally finitely generated} if there exists an open cover
$(U_i)_{i}$ of $\cE$ such that the restriction of $\cE$ to each
$U_i$ is finitely generated.
\end{enumerate}

\subsection{Singular foliations}

As usual (\cf \cite{Stefan, Sussmann}),  a singular foliation will
denote here something more precise than just the partition into
leaves.

\subsubsection{Definitions and examples}

\begin{definition}
Let $M$ be a smooth manifold. A \emph{foliation} on $M$ is a
locally finitely generated submodule of $C_c^\infty(M;TM)$ stable
under Lie brackets.
\end {definition}

\begin{definition}
Let $M,\cF$ be a foliation and $x\in M$. The \emph{tangent space
of the leaf} is the image $F_x$ of $\cF$ in $T_xM$. Put
$I_x=\{f\in C^\infty (M);\ f(x)=0\}$. The \emph{fiber of $\cF$} is
the quotient $\cF_x=\cF/I_x\cF$.
\end {definition}

The evaluation $\tilde e_x:\cF\to T_xM$ vanishes on $I_x\cF$. We
therefore get a surjective homomorphism $e_x:\cF_x\to F_x$ by $e_x
(\xi + f\eta) = \xi(x)$.

The kernel of $\tilde e_x$ is a Lie subalgebra of $\cF$ and
$I_x\cF$ is an ideal in this subalgebra.  It follows that $\ker
e_x=\ker \tilde e_x/I_x\cF$ is a  Lie algebra $\gog_x$.

\begin{examples}\label{exafol}\begin{enumerate}
\item \label{algebroid}Recall that a \emph{Lie algebroid} on $M$
is given by a vector bundle $A$, a Lie bracket $[\ ,\ ]$ on
$C^\infty(M;A)$ together with an \emph{anchor} which is a bundle
map $\#:A\to TM$ satisfying the compatibility relations
$[\#\xi,\#\eta]=\#[\xi,\eta]$ and
$[\xi,f\eta]=f[\xi,\eta]+(\#\xi)(f)\eta$ for all $\xi,\eta\in
C^\infty(M;A)$ and $f\in C^\infty(M)$. To a Lie algebroid there is
naturally associated a foliation: the image of the anchor
$\cF=\#(C_c^\infty(M;A))$. In particular to any Lie groupoid is
associated a foliation.

Note that the anchor maps $I_{x} C_{c}(A)$ onto $I_{x}\cF$, and
since the quotient $C_{c}(A)/ I_{x} C_{c}(A)$ is naturally
identified with the fiber $A_{x}$ we get an onto linear map $A_{x}
\to \cF_{x}$ and the following commutative diagram:
$$
\begin{array}{ccc}
A_{x} &  &   \\
\downarrow & \searrow & \\
\cF_{x} & \stackrel{ev}{\rightarrow} & F_{x}
\end{array}
$$

\item \label{regfol} Recall that a \emph{regular foliation} is a
subbundle $F$ of $TM$ whose sections form a Lie algebroid, \ie
$C^\infty(M;F)$ is stable under Lie brackets. The set of sections
$C_c^\infty(M;F)$ is a foliation in the above sense. We have
$\cF_x=F_x$ for every $x\in M$.

Conversely, Let us show that  if $\cF$ is a foliation such that $F$ is a vector
bundle then $\cF=C_c^\infty (M;F)$. By definition of $F$, we have
$\cF\subset C_c^\infty (M;F)$. Let $x\in M$ and let $\xi_1,\ldots
,\xi_k$ be a basis of $F_x$. There exist $X_1,\ldots ,X_k\in \cF$
such that $X_i(x)=\xi_i$. It follows that $X_1,\ldots ,X_k$ form a
basis of sections of $F$ in a neighborhood of $x$. By a
compactness argument $C_c^\infty (M;F)\subset \cF$.

In particular, in the regular case the distribution $F$ (whence
the partition into leaves) determines the foliation. This is no
longer the case for singular foliations as the following examples
show. \item  Consider the partition of $\R$ into three leaves:
$\R_-^*$, $\{0\}$ and $\R_+^*$. It corresponds to various
foliations $\cF_k$ with $k>0$, where $\cF_k$ is the module
generated by the vector field $x^k\partial/\partial x$. The
foliations $\cF_k$ are different for all different $k$.

In this example, one may legitimately consider that $\cF_1$ is in some sense the best choice of module defining the partition into leaves.

In some other cases there may be no preferred choice of a foliation
(\ie a module) defining the same partition into leaves as the following examples show:

%Let $e_{1},\ldots,e_{n}$ be a base of the Lie algebra $\gog$ and $X_{1},\ldots,X_{n}$ the vector fields in $End(E)$ induced by the infinitesimal action. If $f_{1},\ldots,f_{n}$ are functions on $E$ that vanish at $0$ then $X = \sum f_{i}X_{i}$ vanishes since $\sum f_{i}(0)e_{i} = 0$. This shows that $I_{0}\cF$ vanishes, so $\cF_{0} = \cF$. Let $J_{1}(E)$ be the $1$-jet bundle of $E$, naturally identified with $End(E)$.
%\end{proof}

\item Consider the partition of $\R$ where the
leaves are $\R_+^{*}$ and $\{ x \}$ for every $x \leq 0$. These
are the integral curves of any vector field
$f\frac{\partial}{\partial x}$ where $f(x)$ vanishes for every $x
\leq 0$. So any module generated by such a vector field defines
this foliation. This forces us to impose an a priori choice of
such a module. Note that the module of all vector fields vanishing
in $\R_-$ is not a foliation, as it is not finitely generated.

\item Consider the partition of $\R^2$ into two leaves: $\{0\}$
and $\R^2\setminus \{0\}$. It is given by the action of a Lie
group $G$,  where $G$ can be $GL(2,\R)$, $SL(2,\R)$ or $\C^*$. The
foliation is different for these three actions. The corresponding
$\cF_x$ are equal to $T_x\R^2$ at each non zero $x\in \R^2$, but
$\cF_0$ is the Lie algebra $\gog$. To see this we give the
following more general result:
\end{enumerate}
\end{examples}

\begin{proposition}
Let $E$ be a (real) vector space and $G$ a Lie subgroup of $GL(E)$
with Lie algebra $\gog$. Consider the (linear) action of $G$ on
$E$, and the foliation $\cF$ on $E$ associated with the Lie
groupoid $E\rtimes G$, with Lie algebroid $E\rtimes \gog$. Then
$\cF_{0} = \gog$, \ie the map $p:A_0=\gog\to \cF_{0}$ described in example 1
is injective.
\end{proposition}

\begin{proof}
Let $\tilde \cF\subset C^\infty_c(E;TE)$ be the set of vector
fields on $E$ vanishing at $0$, and $J^1:\tilde \cF\to \cL(E)$ be
the map which associates to a vector field its $1$-jet, \ie its
derivative at $0$. This restricts to $I_{0}\cF\subset I_{0}\tilde
\cF$. It therefore induces a map $\cF_{0} \to \cL(E)$. On the
other hand, $J^1 \circ p:\gog \to \cL(E)$ is the inclusion of
$\gog$ into the Lie algebra of $GL(E)$. It follows that $p$ is injective.
\end{proof}

Let us gather in the following proposition a few easy (and well
known) facts about the fibers and tangent spaces of the leaves.

\begin{proposition} Let $M,\cF$ be a foliation and $x\in M$.\begin{enumerate}
\renewcommand{\theenumi}{\alph{enumi}}
\renewcommand{\labelenumi}{\theenumi)}
\item Let  $X_1,\ldots,X_k\in \cF$ whose images in $\cF_x$ form a
basis of $\cF_x$. There is a neighborhood $U$ of $x$ in $M$ such
that $\cF$ restricted to $U$ is generated by $X_1,\ldots,X_k$.
\item The dimension of $F_x$ is lower semi-continuous and the
dimension of $\cF_x$ is upper semi-continuous. \item The set
$U=\{x\in M;\ e_x:\cF_x\to F_x$ \emph{is bijective}$\}$ is the set
of continuity of $x\mapsto \dim F_x$. It is open and dense in $M$.
The restriction of $F$ to $U$ is a subbundle of $TU$; the
restriction of the foliation $\cF$ to $U$ is the set of sections
of this subbundle. It is a regular foliation.
\end{enumerate}\end{proposition}
\begin{proof}
Since $\cF$ is locally finitely generated there exist
$Y_1,\ldots,Y_N$ which generate $\cF$ in a neighborhood $V$ of
$x$.
\begin{enumerate}
\renewcommand{\theenumi}{\alph{enumi}}
\renewcommand{\labelenumi}{\theenumi)}
\item  Since the images of $X_1,\ldots,X_k$ form a basis of
$\cF_x$, there exist $a_{\ell,i}\in \C$ for $1\le i\le N$ and
$1\le \ell \le k$ such that $Y_i-\sum_{\ell =1}^k a_{\ell,i}X_\ell
\in I_x\cF$. It follows that there exist functions
$\alpha_{i,j}\in C^\infty(M)$ for $1\le i,j\le N$ such that
$\alpha _{i,j}(x)=0$ and for all $i$ we have  $Y_i-\sum_{\ell=1}^k
a_{\ell,i}X_\ell=\sum_{i=1}^n \alpha_{j,i}Y_j$ in a neighborhood
of $x$. This can be written as $\sum_{j=1}^N\beta
_{j,i}Y_j=\sum_{\ell =1}^k a_{\ell,i}X_\ell$ for all $1\le i\le
N$, where $\beta _{i,j}=-\alpha_{i,j}$ if $i\ne j$ and
$\beta_{i,i}=1-\alpha_{i,i}$

For $y\in M$, let $B_y$ denote the matrix with entries
$\beta_{i,j}(y)$. Since $B_x$ is the identity matrix, for $y$ in a
neighborhood $U$ of $x$, the matrix $B(y)$ is invertible. Write
$(B(y))^{-1}=(\gamma_{i,j}(y))$.  We find on $U$,
$Y_i=\sum_{\ell=1}^k c_{\ell,i} X_\ell$, where $c_{\ell,i}=\sum
a_{\ell,j}\gamma_{j,i}$.

\item For $y\in V$, $F_y$ is the image of $\R^N$ through the map
$\varphi_y:(t_1,\ldots ,t_N)\mapsto \sum t_i Y_i(y)$. Since
$y\mapsto \varphi_y$ is continuous, the rank of $F_{y}$ is lower
semi-continuous.

By (a), it follows that on a neighborhood $U$ of $x$, $\cF$ is
generated by $\dim \cF_x$ elements, and therefore, for $y\in U$,
$\dim\cF_y\le \dim \cF_x$. In other words, the dimension of
$\cF_x$ is upper semi-continuous.

\item The function $y\mapsto \dim \cF_y-\dim F_y$ is nonnegative
(since $e_y$ is surjective) and upper semi-continuous by (b). It
follows that $U=\{y;  \dim \cF_y-\dim F_y<1\}$ is open.

Let $V$ be the (open) set where $\dim F_x$ is continuous, \ie has
a local minimum. Note that on $U$ the functions $x\mapsto F_x$ and
$x\mapsto \cF_x$ coincide and are therefore continuous. Therefore
$U\subset V$. Now $F$ restricted to $V$ is a vector bundle, and by
example \ref{exafol}.\ref{regfol}, $\cF$ restricted to $V$ is
equal to $C_c^\infty (V;F)$, whence for every $y\in V$ we have
$\cF_y=F_y$, thus $U=V$.
\end{enumerate}
\end{proof}

\subsubsection{Groups of diffeomorphisms}

Let $M,\cF$ be a foliation. Let $g:M\to N$ be a diffeomorphism.
Associated to $g$ is an isomorphism of modules
$g_*:C_c^\infty(M;TM)\to C_c^\infty(N;TN)$. The image of
$g_*(\cF)$ of $\cF$ is a foliation of $N$. Let us denote by $\cF'$
this foliation. We obviously have  $F'_{g(x)}=dg_x(F_x)$. Also
$g_*(I_x\cF_x)=I_{g(x)}\cF'$, therefore $\cF'_{g(x)}\simeq \cF_x$.

There are several groups of diffeomorphisms of $M$ to be considered:
\begin{enumerate}
\item  The group $\Aut(M,\cF)$ of diffeomorphisms of $M$
preserving the foliation, \ie those diffeomorphisms $g$ such that
$g_*(\cF)=\cF$. \item The group $\exp \cF$ generated by $\exp X$
with $X$ in $\cF$.
\end{enumerate}

The following  fact is fundamental. This is in some sense the main ingredient in the Frobenius theorem:

\begin{proposition} \label{Frob}
We have $\exp \cF\subset \Aut(M,\cF)$. It is actually a normal subgroup in $\Aut(M,\cF)$.
\end{proposition}
\begin{proof}
Let $X\in \cF$; we have to show that $\exp X\in \Aut(M,\cF)$.
Replacing $M$ by a neighborhood of the support of $X$, we may
assume that $\cF$ is finitely generated. Take $Y_1,\ldots ,Y_n$ to
be global sections of $\cF$ generating $\cF$. Since $[X,Y_i]\in
\cF$, there exist functions $\alpha_{i,j}\in C_c^\infty(M)$ such
that $[X,Y_i]=\sum_j \alpha_{j,i}Y_j$.

Denote by $L$ the linear mapping of $C^\infty(M)^n$ given by
$L(f_1,\ldots ,f_n)=(g_1,\ldots ,g_n)$, where
$g_i=X(f_i)+\sum_{j}\alpha_{i,j}f_j$.

Let $S: C^\infty(M)^n\to \cF$ be the map $(f_1,\ldots ,f_n)\mapsto
\sum f_iY_i$; since $L_X\circ S=S\circ L$, we find $\exp X\circ
S=S\circ \exp L$. Therefore, $\exp X(\cF)$, which is the image of
$\exp X\circ S$, is contained in the image of $S$, \ie it is
contained in $\cF$.

Furthermore, if $g\in \Aut (M,\cF)$, we find $g\circ \exp X\circ
g^{-1}=\exp (g_*X)\in \exp \cF$.
\end{proof}

 \begin{definition}
 The \emph{leaves} are the orbits of the action of the group $\exp \cF$.
\end{definition}

\begin{remarks}
\begin{enumerate}
\item For every $x,y$ in the same leaf, there exists (by
definition of the leaves) a diffeomorphism $g\in \exp \cF$ such
that $g(x)=y$. Since $g\in \Aut(M,\cF)$ it follows that the
dimension of $F_x$ and of $\cF_x$ is constant along the leaves.

\item Another group of diffeomorphisms of $\cF$ is the subgroup of
$\Aut(M,\cF)$ consisting of those $g\in \Aut(M,\cF)$ preserving
the leaves, \ie such that $g(x)$ and $x$ are in the same leaf for
all $x\in M$.
\end{enumerate}
\end{remarks}

It will also be useful to consider \emph{local diffeomorphisms}
preserving the foliation. Let $U,V$ be open subsets in $M$. A
local diffeomorphism $\varphi:U\to V$ is said to preserve the
foliation if $\varphi_*(\cF_U)=\cF_V$, \ie if the image through
the diffeomorphism $\varphi$ of the restriction of $\cF$ to $U$ is
the restriction of $\cF$ to $V$. The local diffeomorphisms
preserving the foliation form a pseudogroup. Associated to it is
the \emph{{\'e}tale groupoid of germs of local diffeomorphisms
preserving the foliation}.

\subsubsection{Transverse maps and pull back foliations}

\begin{definition} \label{defip*}
Let $M,N$ be two manifolds, $\varphi :N\to M$ be a smooth map and
$\cF$ be a foliation on $M$.\begin{enumerate}
\renewcommand{\theenumi}{\alph{enumi}}
\renewcommand{\labelenumi}{\theenumi)}
\item Denote by $\varphi ^{-1}(\cF)$ the submodule $\varphi
^{-1}(\cF)=\{X\in C_c^\infty(N;TN)\, ;\ d\varphi (X)\in \varphi
^*(\cF)\}$ of $C_c^\infty(N;TN)$. \item We say that $\varphi $ is
\emph{transverse} to $\cF$ if the natural map $\varphi
^*(\cF)\oplus C_c^\infty(N;TN)\to C_c^\infty(N;\varphi^*(TM))$
(defined by $(\xi,\eta)\mapsto \xi+d\varphi (\eta)$) is onto.
\end{enumerate}
\end{definition}

Of course, a submersion is transverse to any foliation.

\begin{proposition}
Let $M,N$ be two manifolds, $\varphi :N\to M$ be a smooth map and
$\cF$ be a foliation on $M$.
\begin{enumerate}\renewcommand{\theenumi}{\alph{enumi}}
\renewcommand{\labelenumi}{\theenumi)}
\item The $C^\infty(N)$-module $\varphi ^{-1}(\cF)$ is stable
under Lie brackets. \item If $\varphi $ is transverse to $\cF$,
the $C^\infty(N)$-module $\varphi ^{-1}(\cF)$ is locally finitely
generated. It is a foliation.
\end{enumerate}
\end{proposition}
\begin{proof} \begin{enumerate}\renewcommand{\theenumi}{\alph{enumi}}
\renewcommand{\labelenumi}{\theenumi)}
\item Take $X,X'\in f^{-1}(\cF)$ and write $dp(X)=\sum f_iY_i\circ
p$ and $dp(X')=\sum f_i'\,Y_i'\circ p$. We have $dp([X,X'])=\sum
f_if'_j\,[Y_i,Y'_j]\circ p+\sum X(f'_j) Y'_j\circ p-\sum X'(f_i)
Y_i\circ p$. \item By restricting to small open subsets of $N$ and
$M$, we may assume that the tangent bundles of $M$ and $N$ are
trivial and $\cF$ is finitely generated. By projectiveness, there
is a section $s:C_c^\infty(N;p^*(TM))\to p^*(\cF)\oplus
C_c^\infty(N;TN)$. Then, the fibered product
$p^{-1}(\cF)=p^*(\cF)\times_{C_c^\infty(N;p^*(TM))}
C_c^\infty(N;TN)$ identifies with the quotient
$\Big(p^*(\cF)\oplus
C_c^\infty(N;TN)\Big)/s\Big(C_c^\infty(N;p^*(TM))\Big)$ and is
therefore finitely generated.
\end{enumerate}
\end{proof}

Let us note the following obvious facts about pull-back foliations:

\begin{proposition} \label{trivialprop}
Let $M,N$ be two manifolds, $\cF$ be a foliation on $M$ and
$\varphi :N\to M$ be a smooth map transverse to $\cF$. Denote by
$\cF_N$ the pull back foliation.
\begin{enumerate}\renewcommand{\theenumi}{\alph{enumi}}
\renewcommand{\labelenumi}{\theenumi)}
\item For all $x\in N$, we have $(\cF_N)_x=\{\xi \in T_xN;\ d\varphi
_x(\xi)\in F_x\}$. \item Let $P$ be a manifold and $\psi:P\to N$ be
a smooth map. The map $\psi$ is transverse to $\cF_N$ if and only
if $\varphi\circ\psi$ is transverse to $M$ and we have $
(\varphi\circ\psi)^{-1}(\cF)=\psi^{-1}(\cF_N)$.
\end{enumerate}
\ \hfill$\square$
\end{proposition}

\subsection{The leaves and the longitudinal smooth structure}

We now recall the manifold structure of the leaves.

\subsubsection{The local picture}

Let us recall the following well-known fact:

\begin{proposition} \label{localpicture}
Let $(M,\cF)$ be a foliated manifold and let $x\in M$. Put $k=\dim F_x$ and $q=\dim T_xM-\dim F_x$.
\begin{enumerate}\renewcommand{\theenumi}{\alph{enumi}}
\renewcommand{\labelenumi}{\theenumi)}
\item There exist  an open neighborhood $W$ of $x$ in $M$, a foliated manifold $(V,\cF_V)$ of dimension $q$  and a submersion $\varphi:W\to V$ with connected fibers such that $\cF_W=\varphi^{-1}(\cF_V)$ where $\cF_W$ is the restriction of
$\cF$ to $W$. 
\item Moreover,   the tangent space of the leaf of
$(V,\cF_V)$ at the point $\varphi(x)$ is $0$, we have $\ker
(d\varphi)_x=F_x$ and  each fiber of $\varphi$ is contained in a
leaf of $(M,\cF)$.
\end{enumerate}\end{proposition}
\begin{proof}
\begin{enumerate}\renewcommand{\theenumi}{\alph{enumi}}\renewcommand{\labelenumi}{\theenumi)}
\item We proceed by induction in $k$ noting that if $k=0$, there
is nothing to be proved: just take $\varphi:M\to M$ to be the
identity.

Let $X\in \cF$ be a vector field such that $X(x)\ne 0$. Denote by
$\varphi_t=\exp(tX)$ the corresponding one-parameter group. Let
$V_0$ be a locally closed submanifold of $M$ containing $x$ and
such that $T_xV_0\oplus \R X(x)=T_xM$

The map $h:\R\times V_0\to M$ given by $h(t,v)=\varphi_t(v)$ is
smooth and we have $h(0,x)=x$ and $(dh)_0(s,\xi)=sX(x)+\xi$ for
every $s\in \R$ and $\xi\in T_xV_0$. In particular it is bijective,
and it follows that its restriction $\psi:I\times V\to W$ to an
open neighborhood of $(0,x)$ is a diffeomorphism into an open
subset $W$ of $M$. Let $\iota:V\to M$ be the inclusion and
$\varphi:W\to V$ be the composition of $\psi^{-1}:W\to I\times V$
with the projection $I\times V\to V$. Denote by $\cF_W$ the
restriction of $\cF$ to $W$. Note that $V$ is transverse to $\cF$
and put $\cF_V=\iota^{-1}(\cF)$.

Writing $W\simeq I\times V$, every vector field $Y\in
C_c^\infty(W;TW)$ decomposes uniquely as $Y=fX+Z$, where $Z$ is in
the $V$ direction, \ie $Z(s,v)=\varphi_s(Z'_{s}(v))$ with $Z'_s\in
C_c^\infty (V;TV)$. Since $fX\in \cF_W$ and $\varphi_s\in
Aut(\cF)$, it follows that $Y\in \cF_W$ if and only if $Z'_s\in
\cF_V$ for all $s\in I$; in other words, $\cF_W=\varphi^{-1}(V)$.
Now, $(F_V)_x=F_x\cap T_xV$ has dimension $k-1$. We may now apply
the induction hypothesis to $V$.

\item follows from prop. \ref{trivialprop}. To see the last
statement, note that every vector field $Y$ on $W$ tangent to the
fibers of $\varphi :W\to V$ satisfies $d\varphi (Y)=0$, whence
$Y\in \cF_W=\varphi^{-1}(\cF_V)$.
\end{enumerate}
\end{proof}

\subsubsection{Leaves and the smooth structure along the leaves}

Let $M,\cF$ be a foliation. We now describe the topology and smooth structure of a leaf.  We actually consider the collection of leaves as a single manifold with underlying set $M$ \ie a smooth structure on $M$ for which the leaves are open smooth manifolds (of different dimensions).

\begin{definition} Let $N$ be a manifold and $f:N\to M$ be a smooth map. \begin{itemize}
\item We say that $f$ is \emph{leafwise} if $f^{-1}(\cF)=C_c^\infty(N;TN)$ \ie if $df(C_c^\infty(N;TN))\subset
f^*(\cF)$. \item  We say that $f$ is \emph{longitudinally
{\'e}tale} if it is leafwise and for all $x\in N$, the
differential $(df)_x$ is an isomorphism $T_xN\to F_{f(x)}$. \item
A \emph{longitudinal chart} is a locally closed submanifold $U$ of
$M$ such that the inclusion $i:U\to M$ is longitudinally
{\'e}tale.
\end{itemize}
\end{definition}

\begin{proposition} \begin{enumerate}
\renewcommand{\theenumi}{\alph{enumi}}
\renewcommand{\labelenumi}{\theenumi)}
\item There is a new smooth structure on $M$ called the
\emph{leafwise structure}, such that a map $f:N\to M$ (where $N$
is a manifold) is smooth (\resp {\'e}tale) for this structure if
and only if it is smooth and leafwise (\resp longitudinally
{\'e}tale). \item The leaves are the connected components of this
structure.
\end{enumerate}\end{proposition}

\begin{proof}\begin{enumerate}
\renewcommand{\theenumi}{\alph{enumi}}
\renewcommand{\labelenumi}{\theenumi)}

\item We prove that longitudinal charts form an atlas of $M$.

Let $x\in M$. Let $W,V,\varphi $ be as in proposition
\ref{localpicture}. Put $U=\varphi ^{-1}(\{\varphi (x)\})$.

We show that:
\begin{enumerate}\renewcommand{\theenumii}{\roman{enumii}}
\renewcommand{\labelenumii}{(\theenumii)}
\item $U$ is a longitudinal chart;
\item given a leafwise smooth map $g:N\to M$, the set $g^{-1}(U)$ is open in $N$.
\end{enumerate}

Since $\cF_W=\varphi ^{-1}(\cF_V)$, it follows that $g:N\to W$ is
leafwise if and only if $\varphi \circ g$ is. Since the
restriction of $\varphi $ to $U$ is constant, it is leafwise
whence the inclusion $U\to M$ is leafwise. Furthermore, for each
$y\in U$, since it is in the same leaf as $x$, it follows that
$\dim F_y=\dim F_x$, and since $T_yU\subset F_y$, the inclusion
$U\to M$ is longitudinally {\'e}tale. This establishes (i).

Replacing $N$ by its open subset $g^{-1}(W)$, we may assume that
$g(N)\subset W$. We may now replace $W$ by $V$ and $g$ by $\varphi
\circ g$. Therefore (ii) is equivalent to proving: assume $F_x=0$
and $g:N\to M$ is leafwise; then $g^{-1}(\{x\})$ is open. We may
further assume that $N$ is connected and $g^{-1}(\{x\})\ne
\emptyset$. We then prove $g$ is constant.

Let then $y,z\in N$ with $g(y)=x$. Let $y_t$ be a smooth path in
$N$ joining $y$ to $z$. By definition, $(dg)_{y_t}(y'(t))\in
F_{y(t)}$; whence there exists a smooth path $X_t\in \cF$ such
that $X_t(y_t)=(dg)_{y_t}(y'(t))$. Now $x_t=g(y_t)$ as well as the
constant path satisfy a differential equation $x'_t=X_t(x_t)$, and
by uniqueness in Cauchy-Lipschitz theorem they are equal.

Let us now show that longitudinal charts form an atlas. By (i)
they cover $M$. If $U_1$ and $U_2$ are connected longitudinal
charts and $x\in U_1\cap U_2$, let $U$ be as above. It follows by
(ii) that $U\cap U_i$ is open in $U_i$, but since they are
submanifolds with the same dimension, $U\cap U_i$ is open in $U$,
whence $U_1\cap U_2\cap U$ is a neighborhood of $x$ in $U_i$. It
follows that $U_1\cap U_2$ is open in both $U_1$ and $U_2$. This
shows that the longitudinal charts form an atlas.

Now it is quite easy to characterize smooth maps for this new structure.

\item We denote by $N$ the space $M$ endowed with the new structure.

The identity $N\to M$ is leafwise. Let $x\in M$, let $W,V,\varphi
$ be as in proposition \ref{localpicture}, and put $U=\varphi
^{-1}(\{\varphi (x)\})$. By (a), $U$ is open in $N$, and by
proposition \ref{localpicture} it is contained in the leaf $L_x$
of $x$. In particular  $L_x$ is a neighborhood of $x$ in $N$. This
means that the leaves are open in $N$. Therefore, they contain the
connected components of $N$.

On the other hand, let $X_1,\ldots ,X_m\in \cF$ and $x\in M$. The
map $\R\to M$ defined by $f(t)=t\mapsto (\exp tX_1)\circ (\exp
tX_2)\circ\ldots\circ (\exp tX_m)(x)$ is leafwise, therefore $f$
is continuous with values in $N$. It follows that $\exp \cF$ fixes
the connected components of $N$, \ie the leaves are connected.
\end{enumerate}\end{proof}

\begin{remarks} 
\begin{enumerate}
\item Note that from the proof above it follows that if $f:N\to M$
satisfies $(df)_x(T_xN)\subset F_{f(x)}$, then it is leafwise.
\item Let $x\in M$. The tangent space at $x$ for the leafwise
structure is $F_x$. In other words, $F_x$ is the tangent space to
the leaf through $x$.
\end{enumerate}
\end{remarks}

\begin{remark} \label{Algebroidrem}
The foliation $\cF$, when restricted to a leaf gives rise to a transitive algebroid. Indeed, when $M$ is endowed with the leaf topology $(\cF_x)_{x\in M}$ is a vector bundle and the Lie bracket on the sections of this bundle are well defined. More precisely, (and locally) take $x\in M$ and let $V,\cF_V,W$ and $\varphi $ be as in prop. \ref{localpicture}. Put $P=\varphi^{-1}(\{\varphi(x)\})$. Let $\cF_W=C_c^\infty(W)\cF$ be the restriction of $\cF$ to $W$ and put $A_P=\cF_W/I_P\cF_W$ where $I_P=\{f\in C_c(W);\ f_{|P}=0\}$. It is easily seen to be a Lie algebroid over  $P$ which is a neighborhood of $x$ for the leaf topology.
\end{remark}

\section{Bi-submersions and bi-sections}

Throughout this section, we fix a foliation $(M,\cF)$.

\subsection{Bi-submersions}

\begin{definition}  Let $(N,\cF_N)$ be another foliated manifold.
\begin{enumerate}
\renewcommand{\theenumi}{\alph{enumi}}
\renewcommand{\labelenumi}{\theenumi)}\item A
\emph{bi-submersion} between $(M,\cF)$  and $(N,\cF_N)$ is a
smooth manifold $U$ endowed with two smooth maps $s:U\to M$
$t:U\to N$ which are submersions and
satisfying:\begin{enumerate}\renewcommand{\theenumii}{\roman{enumii}}
\renewcommand{\labelenumii}{(\theenumii)}
\item $s^{-1}(\cF)=t^{-1}(\cF_N)$.
\item $s^{-1}(\cF)=C_c^\infty(U;\ker ds)+C_c^\infty (U;\ker dt)$.
\end{enumerate}
\item A \emph{morphism} (\resp a \emph{local morphism}) of
bi-submersions  $(U,t_U,s_U)$ and $(V,t_V,s_V)$ (between $(M,\cF)$
and $(N,\cF_N)$) is a smooth mapping $f:U\to V$ (\resp $f:U'\to V$
where $U'$ is an open subset of $U$) such that for all $u\in U$
(\resp $u\in U'$) we have $s_V(f(u))=s_U(u)$ and
$t_V(f(u))=t_U(u)$. \item A bi-submersion of $(M;\cF)$ is a
bi-submersion between $(M;\cF)$ and $(M;\cF)$. \item A
bi-submersion $(U,t,s)$ of $(M;\cF)$ is said to be
\emph{leave-preserving} if for every $u\in U$, $s(u)$ and $t(u)$
are in the same leaf.
\end{enumerate}\end{definition}

The next proposition motivates the definition of bi-submersions.

\begin{proposition}\label{groupoidbisub}
Let $G$ be a Lie groupoid over $M$ and assume that $\cF$ is the
foliation on $M$ defined by the Lie algebroid $AG$ (example
\ref{exafol}.\ref{algebroid}).  Then $(G,t,s)$ is a bi-submersion
of $(M,\cF)$.
\end{proposition}
\begin{proof}
The Lie algebroid is the restriction $\ker ds|_M$ to $M=G^{(0)}$
of the bundle $\ker ds$. Recall that $\cF=dt(C_c^{\infty}(M;A))$.
In the same way, $\cF=ds(C_c^{\infty}(M;\ker dt|_M))$.

We just have to show that $t^{-1}(\cF)=C_c^\infty(G;\ker
ds)+C_c^\infty (G;\ker dt)$. Using the isomorphism $x\mapsto
x^{-1}$, the equality $s^{-1}(\cF)=C_c^\infty(G;\ker
ds)+C_c^\infty (G;\ker dt)$ will follow.

For $a\in M$, we denote by $G_a=\{x\in G;\ s(x)=a\}$. Let $x\in
G$. The map $y\mapsto yx^{-1}$ is a diffeomorphism between
$G_{s(x)}$ and $G_{t(x)}$, and identifies $(\ker ds)_x$ with
$A_{t(x)}$. In this way, we obtain a bundle isomorphism $\ker
ds\to  t^*(A)$. It follows that $t^*(\cF)=dt(t^*A)=dt (\ker ds)$.
The proposition follows.
\end{proof}

More generally, when the foliation is defined by an algebroid $A$ then any \emph{local Lie groupoid} integrating $A$ gives also a bi-submersion.

Let us come back to the general case. We will use the rather obvious:

\begin{lemma} \label{lembisub} Let  $(U,t,s)$ be a  bi-submersion
of $(M,\cF)$ and  $p:W\to U$ be a submersion. Then $(W,t\circ
p,s\circ p)$ is a bi-submersion.
\end{lemma}
\begin {proof}  The
statement is local and we may therefore assume $W=U\times V$ where
$V$ is a manifold and $p$ is the first projection, in which case
the Lemma is quite obvious.
\end{proof}

\begin{proposition} \textbf{Inverse and composition of bi-submersions.}
Let $(U,t_U,s_U)$ and $(V,t_V,s_V)$ be bi-submersions.
\begin{enumerate}\renewcommand{\theenumi}{\alph{enumi}}
\renewcommand{\labelenumi}{\theenumi)}
\item $(U,s_U,t_U)$ is a bi-submersion. \item Put
$W=U\times_{s_U,t_V}V=\{(u,v)\in U\times V;\ s_U(u)=t_V(v)\}$, and
let $s_W:(u,v)\mapsto s_V(v)$ and $t_W:(u,v)\mapsto t_U(u)$. The
triple $(W,t_W,s_W)$ is a bi-submersion.\end{enumerate}
\end{proposition}
The bi-submersions $(U,s_U,t_U)$ and $(W,t_W,s_W)$ constructed in
this proposition are respectively called the \emph{inverse} of
$(U,t_U,s_U)$ and the \emph{composition} of $(U,t_U,s_U)$ and
$(V,t_V,s_V)$. They are respectively denoted by $U^{-1}$ and
$U\circ V$.

\begin{proof}\begin{enumerate}\renewcommand{\theenumi}{\alph{enumi}}
\renewcommand{\labelenumi}{\theenumi)}
\item is obvious.

\item Denote by $\alpha :W\to M$ the map $(u,v)\mapsto
s_U(u)=t_V(v)$. Note that $p:(u,v)\mapsto u$ and $q:(u,v)\mapsto
v$ are submersions. It follows by lemma \ref{lembisub} that
$(W,t_W,\alpha)$ and $(W,\alpha,s_W)$ are bi-submersions. In
particular $t_W^{-1}(\cF)=\alpha^{-1}(\cF)=s_W^{-1}(\cF)$, whence
$C_c^\infty(W;\ker dt_W)+C_c^\infty(W;\ker ds_W)\subset
t_W^{-1}(\cF)$.
 Note that $\ker d\alpha =\ker dp\oplus \ker dq$. Therefore $$\begin{matrix}
 t_W^{-1}(\cF)&=&C_c^\infty(W;\ker dt_W)+C_c^\infty(W;\ker d\alpha)\hfill&
\\ &=&C_c^\infty(W;\ker dt_W)+C_c^\infty(W;\ker dp)+C_c^\infty(W;\ker dq)& \\
 &=&C_c^\infty(W;\ker dt_W)+C_c^\infty(W;\ker dq)\hfill&\hbox{since $\ker dp\subset \ker dt_W$}\\
 &\subset &C_c^\infty(W;\ker dt_W)+C_c^\infty(W;\ker ds_W)\hfill&
 \end{matrix}$$ since $\ker dp\subset \ker dt_W$. Thus $t_W^{-1}(\cF)=C_c^\infty(W;\ker dt_W)+C_c^\infty(W;\ker ds_W)=s_W^{-1}(\cF)$.
\end{enumerate}
\end{proof}

\begin{remark}
We may actually define a notion of bi-transverse map, more general
than that of bi-submersion. This notion is not used here, but
could be of some help. Let $(N,\cF_N)$ be another foliated
manifold.\begin{enumerate}\renewcommand{\theenumi}{\alph{enumi}}
\renewcommand{\labelenumi}{\theenumi)}
\item A \emph{bi-transverse map} between $(M,\cF)$ and $(N,\cF_N)$
is a smooth manifold $U$ endowed with two smooth maps $s:U\to M$
transverse to $\cF$ and $t:U\to N$ transverse to $\cF_N$
satisfying:\begin{enumerate}\renewcommand{\theenumii}{\roman{enumii}}
\renewcommand{\labelenumii}{(\theenumii)}
\item $s^{-1}(\cF)=t^{-1}(\cF_N)$.
\item $s^{-1}(\cF)=C_c^\infty(U;\ker ds)+C_c^\infty (U;\ker dt)$.
\end{enumerate}
\item A \emph{morphism} (\resp a \emph{local morphism}) of
bi-transverse maps $(U,t_U,s_U)$ and $(V,t_V,s_V)$ is a smooth
mapping $f:U\to V$ (\resp $f:U'\to V$ where $U'$ is an open subset
of $U$) such that for all $u\in U$ (\resp $u\in U'$) we have
$s_V(f(u))=s_U(u)$ and $t_V(f(u))=t_U(u)$.
\end{enumerate}\end{remark}

\subsection{Bisections}

The notion of bisections will allow us to analyze bi-submersions.

\begin{definition} Let  $(U,t,s)$ be a bi-submersion of a $(M,\cF)$.
\begin{enumerate}\renewcommand{\theenumi}{\alph{enumi}}
\renewcommand{\labelenumi}{\theenumi)}
\item A \emph{bisection} of $(U,t,s)$ is a locally closed
submanifold $V$ of $U$ such that the restrictions of both $s$ and
$t$ to $V$ are diffeomorphisms from $V$ onto open subsets of $M$.
\item The \emph{local diffeomorphism associated to a bisection}
$V$ is $\varphi_V =t_V\circ s_V^{-1}$ of $M$ where $s_V:V\to s(V)$
and $t_V:V\to t(V)$ are the restrictions of $s$ and $t$ to $V$.
\item Let $u\in U$ and $\varphi$ a local diffeomorphism of $M$. We
will say that $\varphi$ is carried by $(U,t,s)$ at $u$ if there
exists a bisection $V$ such that $u\in V$ and whose associated
local diffeomorphism coincides with $\varphi$ in a neighborhood of
$u$.
\end{enumerate}
\end{definition}

\begin{proposition}
Let  $(U,t,s)$ be a bi-submersion of a foliation $(M,\cF)$ and
$u\in U$. There exists a bisection $V$ of $(U,t,s)$ containing
$u$.
\end{proposition}

\begin{proof}
The subspaces $\ker (ds)_u$ and $\ker (dt)_u$ of $T_uM$ have the
same dimension. Take a locally closed submanifold $V_0$ of $M$
such that $T_uV_0$ supplements both $\ker (ds)_u$ and $\ker
(dt)_u$ in $T_uM$. One can then take $V$ to be a small
neighborhood of $u$ in $V_0$ - using the local diffeomorphism
theorem.
\end{proof}

The following results are rather obvious. We omit their proofs.

\begin{proposition}
Let $(U,t_U,s_U)$ and $(V,t_V,s_V)$ be
bi-submersions.\label{compbisubm}
\begin{enumerate}\renewcommand{\theenumi}{\alph{enumi}}
\renewcommand{\labelenumi}{\theenumi)}
\item Let $u\in U$ and $\varphi$ a local diffeomorphism of $M$
carried by $(U,t_U,s_U)$ at $u$. Then $\varphi^{-1}$ is carried by
the inverse  $(U,s_U,t_U)$ of $(U,t_U,s_U)$ at $u$. \item Let
$u\in U$ and $v\in V$ be such that $s_U(u)=t_V(v)$ and let
$\varphi,\psi$ a local diffeomorphisms of $M$ carried respectively
by $(U,t_U,s_U)$ at $u$  and by $(V,t_V,s_V)$ at $v$. Then
$\varphi\circ \psi$ is carried by the composition  $(W,t_W,s_W)$
of $(U,t_U,s_U)$ and $(V,t_V,s_V)$ at $(u,v)$.\hfill $\square$
\end{enumerate}
\end{proposition}

\begin{proposition}
If $\varphi:U\to V$ is a local diffeomorphism carried by a
bi-submersion then $\varphi_*(\cF_U)=\cF_V$. \hfill $\square$
\end{proposition}

\subsection{Bi-submersions near the identity; comparison of bi-submersions}

The following result is crucial:

\begin{proposition} \label{identity bisubm}
Let $x\in M$. Let $X_1,\ldots,X_n\in \cF$ be vector fields whose
images form a basis of $\cF_x$. For $y=(y_1,\ldots,y_n)\in \R^n$,
put $\varphi_y=\exp (\sum y_iX_i)\in \exp \cF$. Put $\cW
_0=\R^n\times M$, $s_0(y,x)=x$ and
$t_0(y,x)=\varphi_y(x)$.\begin{enumerate}\renewcommand{\theenumi}{\alph{enumi}}
\renewcommand{\labelenumi}{\theenumi)}
\item There is a neighborhood $\cW $ of $(0,x)$ in $\cW _0$ such
that $(\cW ,t,s)$ is a bi-submersion where $s$ and $t$ are the
restrictions of $s_0$ and $t_0$. \item Let $(V,t_V,s_V)$ be a
bi-submersion and $v\in V$. Assume that $s(v)=x$ and that the
identity of $M$ is carried by $(V,t_V,s_V)$ at $v$. There exists
an open neighborhood $V'$ of $v$ in $V$ and a submersion $g:V'\to
\cW $ which is a morphism of bi-submersions and $g(v)=(0,x)$.
\end{enumerate}
\end{proposition}

\begin{proof} \begin{enumerate}\renewcommand{\theenumi}{\alph{enumi}}
\renewcommand{\labelenumi}{\theenumi)}
\item Consider the vector field $Z:(y,x)\mapsto (0,\sum y_iX_i)$ on $\cW_0$.
Since $Z\in s_0^{-1}(\cF)$, it follows that the diffeomorphism
$\varphi=\exp Z$ fixes the foliation $s_0^{-1}(\cF)$ (prop. \ref{Frob}). Put also
$\alpha:(y,x)\mapsto (-y,x)$ and $\kappa=\alpha\circ \varphi$. These diffeomorphisms also
fix the foliation $s_0^{-1}(\cF)$. Now $\kappa^2=\id$ and $s_0\circ \kappa=t_0$.
It follows  that $s_0^{-1}(\cF)=t_0^{-1}\cF$.

In particular $C_c^\infty(\R^n\times M;\ker dt_0)\subset s_0^{-1}(\cF)$ hence $C_c^\infty(\R^n\times M;\ker ds_0)+C_c^\infty(\R^n\times M;\ker dt_0)\subset s_0^{-1}(\cF)$.

Now, since $\cF$ is spanned by the $X_i$'s near $x$, there exists a smooth function $h=(h_{i,j})$
defined in a neighborhood $W_0$ of $(0,x)$ in $\R^n\times M$ with
values in the space of $n\times n$ matrices such that for $(y,u)\in W$ and $z\in \R^n$ we have:
$(dt_0)_{(y,u)}(z,0)=\sum z_ih_{i,j}(y,u)X_j$, and
$h_{i,j}(0,x)=\delta_{i,j}$. Taking a smaller neighborhood $W$, we
may assume that $h(y,u)$ is invertible. In this neighborhood, $(dt_W)(C_c^\infty(W;\ker ds_W))=t_N^*\cF$, whence
$t_W^{-1}(\cF)\subset C_c^\infty(W;\ker ds_W)+C_c^\infty (W;\ker dt_W)$.

\item Replacing $V$ by an open subset containing $v$, we may
assume that $s_V(V)\subset s(\cW)$ and that the bundles $\ker
dt_V$ and $\ker ds_V$ are trivial. Since $t_V^{-1}(\cF)=
C_c^\infty(V;\ker ds_V)+C_c^\infty (V;\ker dt_V)$, the map
$dt:C_c^\infty(V;\ker ds_V)\to t_V^*(\cF)$ is onto, and there
exist $Y_1,\ldots ,Y_n\in C^\infty(V;\ker ds_V)$ such that
$dt_V(Y_i)=X_i$. Since $X_i(x)$ form a basis of $\cF_x$, the
$Y_i(v)$ are independent. Replacing $V$ by an open neighborhood of
$v$, we may assume that the $Y_i$'s are independent everywhere on
$V$. Let $Z_{n+1},\ldots, Z_k$ be sections of $\ker ds$ such that
$(Y_1,\ldots,Y_n,Z_{n+1},\ldots, Z_k)$ is a basis of $\ker ds_V$.
Since $dt_V(Z_i)\in s_V^*(\cF)$ which is spanned by the $Y_i$'s,
we may subtract a combination of the $Y_i$'s so to obtain a new
basis $(Y_1,\ldots,Y_n,Y_{n+1},\ldots, Y_k)$ satisfying
$t_V^*(Y_i)=X_i$ if $i\le n$ and $t_V^*(Y_i)=0$ if $i>n$. For
$y=(y_1,\ldots,y_k)\in \R^k$ small enough we denote by $\psi_y$
the (partially defined) diffeomorphism $\exp(\sum y_iY_i)$ of $V$.

Let $U_0\subset V$ be a bisection through $v$ representing the
identity, \ie such that  $s_V$ and $t_V$ coincide on $U_0$. We
identify $U_0$ with an open subset of $M$ via  this map. There
exists an open neighborhood $U$ of $v$ in $U_0$ and an open ball
$B$ in $\R^k$ such that $h:(y,u)\mapsto \psi_y(u)$ is a
diffeomorphism of $U\times B$ into an open neighborhood $V'$ of
$v$. Let $p:\R^k\to \R^n$ be the projection to $n$ first
coordinates. The map $p\circ h^{-1}:V'\to \cW$ is the desired
morphism. It is a submersion.
\end{enumerate}
\end{proof}

\begin{corollary} \label{corol}
Let $(U,t_U,s_U)$ and $(V,t_V,s_V)$ be bi-submersions and let
$u\in U$, $v\in V$ be such that
$s_U(u)=s_V(v)$.\begin{enumerate}\renewcommand{\theenumi}{\alph{enumi}}
\renewcommand{\labelenumi}{\theenumi)}
\item If the identity local diffeomorphism is carried by $U$ at
$u$ and by $V$ at $v$, there exists an open neighborhood $U'$ of
$u$ in $U$ and a morphism $f :U'\to V$ such that $f(u)=v$. \item
If there is a local diffeomorphism carried both by $U$ at $u$ and
by $V$ at $v$, there exists an open neighborhood $U'$ of $u$ in
$U$ and a morphism $f :U'\to V$ such that $f(u)=v$. \item If there
is a morphism of bi-submersions $g:V\to U$ such that $g(v)=u$,
there exists an open neighborhood $U'$ of $u$ in $U$ and a
morphism $f :U'\to V$ such that $f(u)=v$.\label{corolc}
\end{enumerate}
\end{corollary}

\begin{proof}
\begin{enumerate}\renewcommand{\theenumi}{\alph{enumi}}
\renewcommand{\labelenumi}{\theenumi)}
\item Put $x=s_U(u)$ and let $\cW $ be as in the proposition
\ref{identity bisubm}. By this proposition there are open
neighborhoods $U'\subset U$ of $u$ and $V'\subset V$ of $v$, and
submersions $g:U'\to \cW$ and $h:V'\to \cW$ which are morphisms
and $g(u)=h(v)=(0,x)$. Let $h_1$ be a local section of $h$ such
that $h'(0,x)=v$. Up to reducing $U'$, we may assume that the
range of $g$ is contained in the domain of $h_1$. Put then
$f=h_1\circ g$. \item Let $\varphi $ be a local diffeomorphism
carried both by $U$ at $u$ and by $V$ at $v$. Up to reducing $U$
and $V$, we may assume that $t_U(U)$ and $t_V(V)$ are contained in
the range of $\varphi $. Then $(U,\varphi^{-1}\circ t_U,s_U)$ and
$(V,\varphi^{-1}\circ t_V,s_V)$ are bi-submersions carrying
respectively at $u$ and $v$ the identity local diffeomorphism. We
may therefore apply (a). \item Let $V_0\subset V$ be a bi-section
through $v$ and $\varphi $ be the local diffeomorphism associated
with $V_0$. Obviously $g(V_0)$ is a bi-section and defines the
same local diffeomorphism. The conclusion follows from (b).
\end{enumerate}\end{proof}

\section{The holonomy groupoid}

\subsection{Atlases and groupoids}

\begin{definition}
Let $\cU=(U_i,t_i,s_i)_{i\in I}$ be a family of bi-submersions.
 \begin{enumerate}\renewcommand{\theenumi}{\alph{enumi}}
\renewcommand{\labelenumi}{\theenumi)}
\item A bi-submersion $(U,t,s)$ is said to be \emph{adapted to
$\cU$ at $u\in U$} if there exists an open subset $U'\subset U$
containing $u$, an element $i\in I$ and a morphism of
bi-submersions $U'\to U_i$. \item A bi-submersion $(U,t,s)$ is
said to be \emph{adapted} to $\cU$ if for all $u\in U$  $(U,t,s)$
is adapted to $\cU$ at $u\in U$. \item We say that
$(U_i,t_i,s_i)_{i\in I}$ is an \emph{atlas} if:\begin{enumerate}
\renewcommand{\theenumii}{\roman{enumii}}
\renewcommand{\labelenumii}{(\theenumii)}
 \item $\bigcup_{i\in I} s_i(U_i)=M$.
 \item The inverse of  every element in $\cU$ is adapted to $\cU$.
 \item The composition of any two elements of $\cU$ is adapted to $\cU$.
 \end{enumerate}
\item Let $\cU=(U_i,t_i,s_i)_{i\in I}$ and
$\cV=(V_j,t_j,s_j)_{j\in J}$ be two atlases. We say that $\cU$ is
\emph{adapted} to $\cV$ if  every element of $\cU$ is adapted to
$\cV$. We say that $\cU$ and $\cV$ are \emph{equivalent} if they
are adapted to each other.
 \end{enumerate}
\end{definition}

\begin{proposition} \textbf{Groupoid of an atlas.}
Let $\cU=(U_i,t_i,s_i)_{i\in I}$ be an atlas.
\begin{enumerate}\renewcommand{\theenumi}{\alph{enumi}}
\renewcommand{\labelenumi}{\theenumi)}
\item The disjoint union $\coprod _{i\in I} U_i$ is endowed with
an equivalence relation $\sim $ given by:  $U_i\ni u\sim v\in U_j$
if there exists a local  morphism from $U_i$ to $U_j$ mapping $u$
to $v$.

Let $G=G_\cU$ denote the quotient of this equivalence relation.
Denote by $Q=(q_i)_{i\in I}:\coprod _{i\in I} U_i\to G$ the
quotient map. \item There are maps $t,s:G\to M$ such that $s\circ
q_i=s_i$ and $t\circ q_i=t_i$. \item For every bi-submersion
$(U,t_U,s_U)$ adapted to $\cU$, there exists a map $q_U:U\to G$
such that, for every local morphism $f:U'\subset U\to U_i$ and
every $u\in U'$, we have $q_U(u)=q_i(f(u))$. \item There is a
groupoid structure on $G$ with set of objects $M$, source and
target maps $s$ and $t$ defined above and such that
$q_i(u)q_j(v)=q_{U_i\circ U_j}(u,v)$.
\end{enumerate}
\end{proposition}

\begin{proof}
It follows from corollary \ref{corol}.\ref{corolc}) that  $\sim $
is an equivalence relation. Assertions (b) and (c) are obvious.

Let $(U,t_U,s_U)$ and $(V,t_V,s_V)$ be bi-submersions adapted to
$\cU$, $u\in U$ and $v\in V$. It follows from corollary
\ref{corol} that $q_U(u)=q_V(v)$ if and only if there exists a
local diffeomorphism carried both by $U$ at $u$ and by $V$ at $v$.
Let $u_i\in U_i$ and $u_j\in U_j$ be such that
$s_i(u_i)=t_j(u_j)$. If $q_U(u)=q_i(u_i)$ and $q_V(v)=q_j(u_j)$,
it follows from proposition \ref{compbisubm} that $(u,v)\in U\circ
V$ and $(u_i,u_j)\in U_i\circ U_j$ carry the same local
diffeomorphism, whence $q_{U\circ V}(u,v)=q_{U_i\circ
U_j}(u_i,u_j)$. Assertion (d) now follows.
\end{proof}

\begin{remarks} \begin{enumerate}
\renewcommand{\theenumi}{\alph{enumi}}
\renewcommand{\labelenumi}{\theenumi)}
\item Let $\cU$ and $\cV$ be two atlases. If $\cU$ is adapted to
$\cV$ there is a natural morphism of groupoids $G_\cU\to G_\cV$.
This homomorphism is injective. It is surjective if and only if
$\cU$ and $\cV$ are equivalent.

\item We may also consider a bi-submersion $(U,t,s)$ defined on a
\emph{not necessarily Hausdorff} manifold $U$. Since every point
in $U$ admits a Hausdorff neighborhood, we still get a map
$q_U:U\to G$. \end{enumerate}\end{remarks}

\begin{examples} \label{examples}
There are several choices of non equivalent atlases that can be constructed.
\begin{enumerate}
\item The \emph{full holonomy atlas}. Take all possible
bi-submersions.\footnote{To make them a set, take bi-submersions
defined on open subsets of $\R^N$ for all $N$.} Of course, any
bi-submersion is adapted to this atlas and the corresponding
groupoid is the biggest possible. \item The \emph{leaf preserving
atlas.}  Take all leave preserving bi-submersions, \ie those
bi-submersions $(U,t,s)$ such that for all $u\in U$, $s(u)$ and
$t(u)$ are in the same leaf. It is immediately seen to be an
atlas. \item Take a cover of $M$ by $s$-connected bi-submersions
given by proposition \ref{identity bisubm}. Let $\cU$ be the atlas
generated by those. This atlas is adapted to any other atlas and
the corresponding groupoid is the smallest possible. We say that
$\cU$ is   a \emph{path holonomy atlas.} The associated maximal
atlas will be called the path holonomy atlas.

\item \label{exa4}Assume that the foliation $\cF$ is associated
with a (not necessarily Hausdorff) Lie groupoid $(G,t,s)$. Then
$(G,t,s)$ is a bi-submersion (prop. \ref{groupoidbisub}). It is
actually an atlas. Indeed, $G^{-1}$ is isomorphic to $G$ via
$x\mapsto x^{-1}$ and $G\circ G=G^{(2)}$ is adapted to $G$ via the
composition map $G^{(2)}\to G$ which is a morphism of
bi-submersions.

If $G$ is $s$-connected, this atlas is equivalent to the path
holonomy atlas. Indeed, consider the Lie algebroid $AG$ associated
to $G$. Let $U$ be an open neighborhood of a point $x \in M$ and
choose a base of sections $\xi_{1},\ldots,\xi_{n}$ of $A_{U}$.
Then the exponential map gives rise to a diffeomorphism $U \times
B^n \to W$, where $B^{n}$ is an open ball at zero in $\R^n$ and
$W$ is an open subset at $1_{x} \in G$. This map is
$(x,\theta_1,\ldots,\theta_n) \mapsto \exp(\sum \theta_i
\xi_{i})(1_x)$. (For an account of the exponential map, as well as
this diffeomorphism see \cite{Nistor}.)

The restrictions of the source and target maps of $G$ make $W$ a
bi-submersion. On the other hand $U \times B^n$ is a bi-submersion
with source the first projection and target the map
$(x,\theta_1,\ldots,\theta_n) \mapsto \exp(\sum \theta_i
\alpha(\xi_i))(x)$, where $\alpha$ is the anchor map. It is then
straightforward to check that the diffeomorphism above is actually
an isomorphism of bi-submersions.

Since $G$ is $s$-connected, it is generated by $W$, and it follows
that the atlas $\{(G,t,s)\}$ is equivalent to the atlas generated
by $(U_i \times B^n,s,t)_{i \in I}$ \ie the path holonomy atlas.
The path holonomy groupoid is therefore a quotient of $G$.
\end{enumerate}
\end{examples}

\begin{definition}
The \emph{holonomy groupoid} is the groupoid associated with the
path holonomy atlas.
\end{definition}

\begin{remark}
Let $x\in M$ and $g\in \exp\cF$. There is a bi-submersion carrying
$g$ at $x$ adapted to the path holonomy atlas. In other words, we
get a surjective morphism of groupoids $\exp \cF\ltimes M\to G$.
If $g\in \exp (I_x\cF)$, then the image of $(g,x)$ is equal to
that of $(\id,x)$.
\end{remark}

\begin{examples} 
\begin{enumerate}
\item Let us illustrate in a simple case how bad the topology of the holonomy groupoid is when the foliation is not almost regular. Consider the foliation of $\R^2$ with two leaves $\{0\}$ and $\R^2\setminus \{0\}$ given by the action of $SL(2,\R)$. Recall that the groupoid $\R^2\rtimes SL(2,\R)$ is, as a set $SL(2,\R)\times \R^2$, that we have $s(g,u)=u$ and $t(g,u)=gu$ for $u\in \R^2$ and $g\in SL(2,\R)$. It is a bi-submersion and an atlas for our foliation; this atlas is equivalent to the path holonomy atlas since it is $s$-connected. The path holonomy groupoid $G$ is therefore a quotient of this groupoid.

Let us prove that $G=(SL(2,\R)\times \{0\})\cup (\R^2\setminus \{0\})\times (\R^2\setminus \{0\})$.  Note that $(\R^2\setminus \{0\})\times (\R^2\setminus \{0\})$ is a bi-submerrsion. We therefore just have to show that the map $(SL(2,\R)\times \{0\})\to G$ is injective.

Let $g\in SL(2,\R)$ and let $V$ be a bi-section through $(g,0)$. In particular, $s:V\to \R^2$ is a local diffeomorphism, so that $V$ is the graph of a smooth map $\varphi :B\to SL(2,\R)$, where $B$ is a neighborhood of $0$ in $\R^2$, such that $\varphi(0)=g$. The partial diffeomorphism defined by $V$ is then $u\mapsto \varphi(u)u$ and its derivative at $0$ is $g$. In particular, two distinct elements $g, g'$ in $SL(2,\R)$ don't give rise to the same bisection. Therefore $(g,0$ and $(g',0)$ are different.

Let now $x\in \R^2\setminus \{0\}$ and $g\in SL(2,\R)$ such that $gx=x$. The sequence $\Big(\frac{x}{n},\frac{x}{n}\Big)\in G$ converges in $G$ to $(g,0)$. In other words, this sequence converges to every point of the stabilizer of $x$: the set of its limits is a whole real line!

\item More generally, let $E$ be a finite dimensional real vector space and $\Gamma$ a closed connected subgroup of $GL(E)$. Consider the group $\Gamma$ acting on the space $E$ and the associated Lie groupoid $E\rtimes \Gamma$. It defines a foliation $\cF$ on $E$. The holonomy groupoid $G$ is again a quotient of this groupoid. As above, one sees that $(g,0)\sim (g',0) \iff g=g'$ (for $g,g'\in \Gamma$).

\end{enumerate}
\end{examples}

\subsection{The ``quasi-regular'' case}

The holonomy groupoid was first defined by Winkelnkemper in the
regular case (\cf \cite{Winkelnkemper:holgpd}) and generalized to singular cases by
various authors.   In particular, a construction suggested by
 Pradines and Bigonnet  \cite{Pradines,
Bigonnet-Pradines} was carefully analyzed, and its
precise range of applicability found, by Claire Debord \cite{Debord:2001LILA, Debord:2001HGSF}. Their construction holds for foliations
which are locally projective \ie when $\cF$ is a Lie algebroid
(\cf \cite{Debord:2001LILA, Debord:2001HGSF}). It  follows that
such a Lie algebroid is integrable (a result given also by Crainic
and Fernandes in \cite{Crainic-Fernandes}). One can summarize
their result in the following way:

\begin{proposition} \label{prop3.9}
(\cf \cite{Winkelnkemper:holgpd, Pradines, Bigonnet,
Bigonnet-Pradines, Debord:2001LILA, Debord:2001HGSF}) Assume that
$\cF$ is a Lie algebroid. Then there exists an $s$-connected Lie
groupoid $G$ with Lie algebroid $\cF$ which is minimal in the
sense that if $H$ is an $s$-connected Lie groupoid  with Lie
algebroid $\cF$, there is a surjective {\'e}tale homomorphism of
groupoids $H\to G$ which is the identity at the Lie algebroid
level. Moreover, $G$ is a quasi-graphoid.\hfill $\square$
\end{proposition}

Recall, \cf \cite{Bigonnet,  Debord:2001HGSF} that a
\emph{quasi-graphoid} is a Lie groupoid $G$ such that if $V\subset
G$ is a locally closed submanifold such that $s$ and $t$ coincide
on $V$ and  $s:V\to M$ is {\'e}tale, then $V\subset G^{(0)}$.

It is now almost immediate that the holonomy groupoid given here
generalizes the holonomy groupoid of Winkelnkemper,
Pradines-Bigonnet and Debord.

\begin{proposition} \label{prop3.10}
Let $G$ be an $s$-connected quasi-graphoid and let $\cF$ be the
associated foliation. Then the holonomy groupoid of $\cF$ is
canonically isomorphic to $G$.
\end{proposition}

\begin{proof}
By proposition \ref{groupoidbisub} and  example
\ref{examples}.\ref{exa4}, it follows that $G$ is a bi-submersion
associated with $\cF$, and an atlas equivalent to the path
holonomy atlas $\cU$ (since $G$ is $s$-connected). We thus get a
surjective morphism $G\to G_\cU$. This morphism is injective by
the definition of quasi-graphoids.
\end{proof}

Putting together prop. \ref{prop3.9} and prop. \ref{prop3.10}, we find:

\begin{corollary}
In the regular case, our groupoid coincides with the one defined
by Winkelnkemper; in the quasi-regular case, it coincides with the
one defined by Pradines-Bigonnet and Debord.  \hfill $\square$
\end{corollary}

\subsection{The longitudinal smooth structure of the holonomy groupoid}

The topology of the holonomy groupoid, as well as all groupoids
associated with other atlases, is usually quite bad. Let us discuss
the smoothness issue.

%On the other hand, it carries a natural \emph{longitudinal} smooth structure.

Fix an atlas $\cU$ and let $G$ be the associated groupoid. For
every bi-submersion $(U,t,s)$ adapted to $\cU$, let
$\Gamma_U\subset U$ be the set of $u\in U$ such that $\dim
T_uU=\dim M+\dim \cF_{s(u)}$. It is an open subset of $U$ when $U$
is endowed with the smooth structure along the leaves of the
foliation $t^{-1}(\cF)=s^{-1}(\cF)$.

\begin{proposition} \label{smoothlong}
\begin{enumerate} \renewcommand{\theenumi}{\alph{enumi}}
\renewcommand{\labelenumi}{\theenumi)}
\item For every $x\in G$, there is a  bi-submersion $(U,t,s)$
adapted to $\cU$ such that $x\in q_U(\Gamma_U)$. \item
\label{itembb} Let $(U,t,s)$ and $(U',t',s')$ be two
bi-submersions and let $f:U\to U'$ be a morphism of
bi-submersions. Let  $u\in U$. If $u\in \Gamma_U$, then $(df)_u$
is injective; if $f(u)\in \Gamma_{U'}$, then $(df)_u$ is
surjective.
\end{enumerate}

\end{proposition}
\begin{proof} \begin{enumerate} \renewcommand{\theenumi}{\alph{enumi}}
\renewcommand{\labelenumi}{\theenumi)}
\item Let $x\in G$. Let $(W,t,s)$ be a bi-submersion adapted to
$\cU$ and $w\in W$ and such that $x=q_W(w)$. Let $A\subset W$ be a
bi-section through $w$. Let $g:s(A)\to t(A)$ be the partial
diffeomorphism of $M$ defined by $A$. By proposition \ref{identity
bisubm} there exists a bi-submersion $(U_0,t_0,s_0)$ and $u_0\in
U_0$ such that $\dim U_0=\dim \cF_{s(u)}+\dim M$, $s_0(u_0)=s(u)$
and carrying the identity through $u_0$. Put then $U=\{x\in U_0;\
t_0(u)\in s(A)\}$ and let $s$ be the restriction of $s_0$ to $U$
and $t$ be the map $u\mapsto g(t_0(u))$. Obviously $(U,t,s)$ is a
bi-submersion which carries $g$ at $u_0$. It follows that
$(U,t,s)$ is adapted to $\cU$ at $u_0$ and $q_U(u_0)=q_W(w)=x$. It
is obvious that $u_0\in \Gamma_U$. \item Since $s$ and $s'$ are
submersions and $s'\circ f=s$, $df_u$ is injective or surjective
if and only if its restriction $(df_{|\ker ds})_u:\ker (ds)_u\to
\ker (ds')_{f(u)}$ is. Consider the composition $$\ker
(ds)_u\buildrel{(df_{|\ker ds})_u}\over {\longrightarrow}\ker
(ds')_{f(u)}\buildrel{t'_*}\over {\longrightarrow}{} \cF_{s(u)}.$$
By definition of bi-submersions the maps $t'_*$ and $t_*=t'_*\circ
(df)_u$ are onto; if $u\in \Gamma_U$, then $t_*:\ker (ds)_u\to
\cF_{s(u)}$ is an isomorphism;  if $f(u)\in \Gamma_{U'}$, then
$t'_*:\ker (ds')_{f(u)}\to \cF_{s(u)}$ is an isomorphism. The
conclusion follows.
\end{enumerate}\end{proof}

It follows from \ref{smoothlong}.\ref{itembb}) that the
restriction of $f$ to a neighborhood of $\Gamma_U\cap
f^{-1}(\Gamma_{U'})$ is {\'e}tale. This restriction preserves the
foliation, and is therefore {\'e}tale also with respect to this
structure. Now the $\Gamma_U$ are open in $U$ with respect to the
longitudinal structure; they are manifolds. The groupoid $G$ is
obtained by gluing them through local diffeomorphisms.

\medskip
Let us record for further use the following immediate consequence
of prop. \ref{smoothlong}.

\begin{corollary} \label{corsmoothlong}
Let $(U,t,s)$ and $(U',t',s')$ be two bi-submersions, let $f:U\to
U'$ be a morphism of bi-submersions. Let  $u\in U$. There exist
neighborhoods $V$ of $u$ in $U$, $V'$ of $f(u)$ in $U'$, a
bi-submersion $(U'',t'',s'')$ and a morphism $p:V'\to U''$ which
is a submersion, such that $f(V)\subset V'$ and the map $v\mapsto
p(f(v))$ is also a submersion.
\end{corollary}

\begin{proof}
Just take a bi-submersion $(U'',t'',s'')$ and $u''\in
\Gamma_{U''}$ such that $q_U(u)=q_{U''}(u'')$. By definition of
the groupoid, there exists a morphism $p$ of bi-submersions from a
neighborhood of $f(u)$ to $U''$, it follows from prop.
\ref{smoothlong}.\ref{itembb}) that up to restricting to small
neighborhoods, one may assume that $p$ and $p\circ f$ are
submersions.
\end{proof}

\begin{remark}
In many cases, the groupoid $G_\cU$ is a manifold when looked at longitudinally. We observe that the necessary and sufficient condition for it to be a manifold is the following: We need to
ensure that for every $x \in M$ there exists a bi-submersion
$(U,t,s)$ in the path holonomy atlas and a $u \in U$ which carries
the identity at $x$ and has an open neighborhood $U_{u} \subseteq
U$ with respect to the leaf topology such that the quotient map
$U_{u} \to G_{\cU}$ is injective. Under this condition the $s$
($t$)-fibers of $G_{\cU}$ are smooth manifolds (of dimension equal
to the dimension of $\cF_x$) and the quotient map $q_{U} : U \to
G_{\cU}$ is a submersion for the leaf smooth structure. Note that this condition does not depend on the atlas $\cU$.

If this condition is satisfied, we say that the groupoid $G_\cU$ is \emph{longitudinally smooth.}

Let us note that if $G_\cU$ is longitudinally smooth, then its restriction to each leaf is a Lie groupoid (with respect to the  smooth structure of the leaf).  The Lie algebroid of this groupoid is then the transitive Lie algebroid considered in remark \ref{Algebroidrem}. So a necessary condition is that this algebroid be integrable. Since this algebroid is transitive, a necessary and sufficient condition is the one given by Crainic and Fernandes in \cite{Crainic-Fernandes} following Mackenzie (\cite{Mackenzie}), namely that a given subgroup of the center of the Lie algebra $\gog_x(=\ker (\cF_x\to F_x)$) be discrete.

We therefore come up with the two following questions:\begin{enumerate}
\item Could it happen that the above Mackenzie-Crainic-Fernandes obstruction always vanishes for the algebroid associated with a foliation?
\item If this obstruction vanishes, does this imply that our groupoid is longitudinally smooth?
\end{enumerate}
\end{remark}

In case the groupoid $G_\cU$ is not longitudinally smooth, we are led to
give the following definition:

\begin{definition}
A {\em holonomy pair} for a foliation $(M, \cF)$ is a pair $(\cU,
G)$ where $\cU$ is an atlas of bi-submersions and $G$ is a
groupoid over $M$ which is a Lie groupoid for the smooth
longitudinal structure, together with a a surjective groupoid
morphism $\alpha:G_\cU\to G$ such that the maps $\alpha \circ
q_{U}$ are leafwise submersions  for each $U\in \cU$.
\end{definition}
We will describe in section 4 how to construct 
reduced $C^{*}$-algebras associated with a given holonomy pair. In
case $(\cU,G_{\cU})$ is not a holonomy pair for some atlas $\cU$
we can always replace $G_{\cU}$ with the groupoid $R_{\cF}$
defined naturally by the equivalence relation of ``belonging in
the same leaf" (or leaves that are related by $\cU$ -  if $\cU$ is
not leaf preserving). This groupoid is not smooth in general, but it is always leafwise
smooth: it is the disjoint union of $L\times L$ where $L$ is a
leaf (or the disjoint union of $\cU$-equivalent leaves). 

This will allow us to deal with the cases where the atlas $\cU$ does not satisfy the
condition we mentioned above.

\section{The $C^*$-algebra of the foliation}

We will construct a convolution $*$-algebra on the holonomy
groupoid $G$. We then construct a family of regular
representations, which will yield the reduced $C^*$-algebra of the
foliation. Moreover, we construct an $L^1$-norm, and therefore a
full $C^*$-algebra of the foliation.

\subsection{Integration along the fibers of a submersion}

Let $E$ be a vector bundle over a space $X$ and $\alpha \in \R$.
We denote by $\Omega ^\alpha E$ the bundle of $\alpha $-densities
on $E$.

Let $\varphi:M\to N$ be a submersion, and let $E$ be a vector
bundle on $N$. Integration along the fibers gives rise to a linear
map $\varphi_!:C_c(M;\Omega^1(\ker d\varphi)\otimes \varphi^*E)\to
C_c(N;E)$ defined by $\varphi_!(f):x\mapsto
\int_{\varphi^{-1}(x)}f\ \ (\in E_x)$.

The following result is elementary:

\begin{lemma} \label{obviouslemma}
With the above notation, if $\varphi :M\to N$ is a surjective
submersion, then the map $\varphi_!:C_c(M;\Omega^1(\ker
d\varphi)\otimes \varphi^*E)\to C_c(N;E)$ is surjective.
\end{lemma}

\begin{proof}
This is obvious if $M=\R^k\times N$. The general case follows
using partitions of the identity.
\end{proof}
\subsection{Half densities on bi-submersions}

Let $(U,t_U,s_U)$ be a bi-submersion. Consider the bundle $\ker
dt_U\times \ker ds_U$ over $U$. Let $\Omega^{1/2}U$ denote the
bundle of (complex) half densities of this bundle.

Recall that the inverse bi-submersion $U^{-1}$ of $U$ is
$(U,s_U,t_U)$. Note that $\Omega^{1/2}U^{-1}$ is canonically
isomorphic to $\Omega^{1/2}U$.

Let $(V,t_V,s_V)$ be another bi-submersion. Recall that the
composition $U\circ V$ is the fibered product $U\times
_{s_U,t_V}V$. Denote by $p:U\circ V\to U$ and $q:U\circ V\to V$
the projections. By definition of the fibered product, we have
identifications: $\ker dp\sim q^*\ker dt_V$ and $\ker dq\sim
p^*\ker ds_U$. Furthermore, since $t_{U\circ V}=t_U\circ p$ and
$s_{U\circ V}=s_V\circ q$, we have exact sequences $$0\to \ker
dp\to \ker dt_{U\circ V}\to p^*(\ker dt_U)\to 0\ \ \ \hbox{and}\ \
\  0\to \ker dq\to \ker ds_{U\circ V}\to q^*(\ker ds_V)\to 0.$$ In
this way, we obtain a canonical isomorphism $\Omega^{1/2}(U\circ
V)\simeq p^*(\Omega^{1/2}U)\otimes q^*(\Omega^{1/2}V)$.

\begin{definition} \label{def4.2}
Let $(U,t_U,s_U)$ and $(V,t_V,s_V)$ be bi-submersions.
\label{prod*}\begin{enumerate}
\renewcommand{\theenumi}{\alph{enumi}}
\renewcommand{\labelenumi}{\theenumi)}
\item Denote by $U^{-1}$ the inverse bi-submersion and $\kappa
:U\to U^{-1}$ the (identity) isomorphism. For $f\in C_c^\infty
(U;\Omega^{1/2}U)$ we set $f^*=\overline f\circ \kappa^{-1}\in
C_c^\infty (U^{-1};\Omega^{1/2}U^{-1})$ - via the identification
$\Omega^{1/2}U^{-1}\sim \Omega^{1/2}U$. \item Denote by $U\circ V$
the composition. For $f\in C_c^\infty (U;\Omega^{1/2}U)$ and $g\in
C_c(V;\Omega^{1/2}V)$ we set $f\otimes g:(u,v)\mapsto f(u)\otimes
g(v) \in C_c^\infty (U\circ V;\Omega^{1/2}U\circ V)$ - via the
identification $\Omega^{1/2}(U\circ V)_{(u,v)}\simeq
(\Omega^{1/2}U)_u\otimes (\Omega^{1/2}V)_v$.
\end{enumerate} \end{definition}

Let $(U,t_U,s_U)$ and $(V,t_V,s_V)$ be bi-submersions and let
$p:U\to V$ be a submersion which is a morphism of bi-submersions.
We have exact sequences $$0\to \ker dp\to \ker ds_U\to p^*(\ker
ds_V)\to 0\ \ \hbox{and}\ \ 0\to \ker dp\to \ker dt_U\to p^*(\ker
dt_V)\to 0$$ of vector bundles over $U$, whence a canonical
isomorphism $$\Omega^{1/2}(U)\simeq (\Omega ^{1/2} \big((\ker dp
\oplus p^*(\ker dt_U )\oplus (\ker dp\oplus  p^*(\ker
ds_U)\big)\simeq (\Omega ^{1}\ker dp)\otimes p^*(\Omega
^{1/2}V).$$ Integration along the fibers then gives a map
$p_!:C_c^\infty (U;\Omega^{1/2}U)\to C_c^\infty
(V;\Omega^{1/2}V)$.

\subsection{The $*$-algebra of an atlas}

Let us now fix an atlas $\cU=(U_i,t_i,s_i)_{i\in I}$. Denote by
$U$ the disjoint union $\coprod_{i\in I} U_i$, and let
$t_U,s_U:U\to M$ be the submersions whose restrictions to $U_i$ are
$t_i$ and $s_i$ respectively. Then $(U,t_U,s_U)$ is a
bi-submersion, and $C_c^\infty (U;\Omega^{1/2}U)=\bigoplus _{i\in
I}C_c^\infty (U_i;\Omega^{1/2}U_i)$.

\begin{lemma} \label{Lemma4.2}
Let $(V,t_V,s_V)$ be a bi-submersion adapted to
$\cU$.\begin{enumerate} \renewcommand{\theenumi}{\alph{enumi}}
\renewcommand{\labelenumi}{\theenumi)}
\item Let $v\in V$. Then there exist a bi-submersion $(W,t_W,s_W)$
and submersions $p:W\to U$ and $q:W\to V$ which are morphisms of
bi-submersions such that $v\in q(W)$.
\item Let $f\in C_c^\infty (V;\Omega^{1/2}V)$. Then there exist a
bi-submersion $(W,t_W,s_W)$, submersions $p:W\to U$ and $q:W\to V$
which are morphisms of bi-submersions and $g\in C_c^\infty
(W;\Omega^{1/2}W)$ such that $q_!(g)=f$.
\end{enumerate}\end{lemma}

\begin{proof}
\begin{enumerate} \renewcommand{\theenumi}{\alph{enumi}}
\renewcommand{\labelenumi}{\theenumi)}
\item Let $h$ be the local diffeomorphism carried by $V$ at $v$. Since $V$ is adapted to $\cU$, there exists $u\in U$ such
that  $h$ is also carried by $U$ at $u$. By prop.
\ref{smoothlong}, there exist neighborhoods $V'$ of $v\in V$ and
$U'$ of $u$ in $U$ a bi-submersion $(N,t_N,s_N)$ and morphisms
$\alpha :U'\to N$ and $\beta :V'\to N$ such that $\alpha (u)=\beta
(v)\in \Gamma_N$. Taking smaller neighborhoods, we may further
assume that $\alpha $ and $\beta $ are submersions. Let $W$ be the
fibered product $U'\times _NV'$.

\item Using (a), we find bi-submersions $(W_j,t_j,s_j)$ and
submersions which are morphisms of bi-submersions $p_j:W_j\to U$
and $q_j:W_j\to V$ such that $\bigcup q_j(W_j)$ contains the
support of $f$. Put $W=\coprod W_j$. It follows from Lemma
\ref{obviouslemma} that $f$ is of the form $q_!(g)$.
\end{enumerate}
\end{proof}

Our algebra  $\cA_\cU$ is the quotient
$\cA_\cU=C_c(U;\Omega^{1/2}U)/\cI$ of
$C_c(U;\Omega^{1/2}U)=\bigoplus _{i\in I}C_c(U_i;\Omega^{1/2}U_i)$
by the subspace $\cI $ spanned by the $p_!(f)$ where $p:W\to U$ is
a submersion and $f\in C_c(W;\Omega^{1/2}W)$ is such that there
exists a morphism $q:W\to V$ of bi-submersions which is a
submersion and such that $q_!(f)=0$.

From lemma \ref{Lemma4.2}, it follows:

\begin{proposition}
To every bi-submersion $V$ adapted to $\cU$ one associates a
linear map $Q_V:C_c(V;\Omega^{1/2}V)\to \cA_\cU$. The maps $Q_V$
are characterized by the following properties:\begin{enumerate}
\renewcommand{\theenumi}{\roman{enumi}}
\renewcommand{\labelenumi}{(\theenumi)}
\item If $(V,t_V,s_V)=(U_i,t_i,s_i)$, $Q_V$ is the quotient map \\
$C_c(U_i,\Omega^{1/2}U_i)\subset \bigoplus _{j\in
I}C_c(U_i;\Omega^{1/2}U_j)\to \cA_\cU=\bigoplus _{i\in
I}C_c(U_i;\Omega^{1/2}U_i)/\cI$. \item For every morphism $p:W\to
V$ of bi-submersions which is a submersion, we have $Q_W=Q_V\circ
p_!$.
\end{enumerate}
\end{proposition}

\begin{proof}
With the notation of lemma \ref{Lemma4.2}.b), we define $Q_V(f)$
to be the class of $p_!(g)$. This map is well defined, since if we
are given two bi-submersions $W$ and $W'$ with morphisms $p:W\to
U$, $q:W\to V$, $p':W'\to U$ and $q':W'\to V$ which are
submersions and elements $g\in C_c(W,\Omega^{1/2}W)$ and $g'\in
C_c(W',\Omega^{1/2}W')$ such that $q_!(g)=q'_!(g')$, we let $W''$
to be the disjoint union of $W$ and $W'$ and $g''$ to agree with
$g$ on $W$ and with $-g'$ in $W'$. With obvious notations, we have
$q''_!(g'')=q_!(g)-q'_!(g')=0$, hence
$p_!(g)-p'_!(g')=p''_!(g'')=0$. It is obvious that this
construction satisfies (i) and (ii) and is characterized by these
properties.
\end {proof}

\begin{proposition}
There is a well defined structure of $*$-algebra on $\cA_\cU$ such
that, if $V$ and $W$ are bi-submersions adapted to $\cU$, $f\in
C_c(V;\Omega^{1/2}V)$ and $g\in C_c(W;\Omega^{1/2}W)$, with the
notation given in definition \ref{def4.2}, we have
$$\Big(Q_{V}(f)\Big)^*= (Q_{V^{-1}})(f^*) \ \ \hbox{and}\ \
Q_{V}(f)Q_{W}(g)=Q_{V\circ W}(f\otimes g).$$
\end{proposition}

\begin{proof}
If $p:V\to V'$ is a submersion which is a morphism of
bi-submersions, and if $f\in C_c(V;\Omega^{1/2}V)$ satisfies
$p_!(f)=0$, then $(p\times \id_W)_!(f\otimes g)=0$, where
$(p\times \id_W):V\circ W\to V'\circ W$ is the map $(v,w)\mapsto
(p(v),w)$. Together with a similar property, given a submersion
$p:W\to W'$ which is a morphism of bi-submersions, this shows that
the product is well defined.
\end{proof}

\begin{definition} \label{*algebra}
The $*$-algebra of the atlas $\cU$ is $\cA_\cU$ endowed with these
operations. The $*$-algebra $\cA(\cF)$ of the foliation is
$\cA_\cU$ in the case $\cU$ is the path holonomy atlas.
\end{definition}

\subsection{The $L^1$-norm and the full $C^*$-algebra}

Let $(U,t,s)$ be a bi-submersion, $u\in U$ and put $x=t(u)$, and
$y=s(u)$. The image by the differential $dt_u$ of $\ker ds$ is
$F_y$ and the image by the differential $ds_u$ of $\ker dt$ is
$F_x$. We therefore have exact sequences of vector spaces $$0\to
(\ker ds)_u\cap (\ker dt)_u\to (\ker ds)_u\to F_y\to 0\
\hbox{and}\ 0\to (\ker ds)_u\cap (\ker dt)_u\to (\ker dt)_u\to
F_x\to 0.$$ We therefore get isomorphisms $$\Omega^{1/2}(\ker
ds)_u\simeq \Omega^{1/2}\Big((\ker ds)_u\cap (\ker
dt)_u\Big)\otimes \Omega^{1/2}F_y\ \hbox{and} \ \Omega^{1/2}(\ker
dt)_u\simeq \Omega^{1/2}\Big((\ker ds)_u\cap (\ker
dt)_u\Big)\otimes \Omega^{1/2}F_x$$

Choose a Riemannian metric on $M$. It gives a Euclidean metric on
$F_x$ and $F_y$ and whence trivializes the one dimensional vector
spaces $\Omega^{1/2}F_x$ and $\Omega^{1/2}F_y$. We thus get an
isomorphism $\rho^U_u:\Omega^{1/2}(\ker ds)_u\to \Omega^{1/2}(\ker
dt)_u$.

Let $f$ be a nonnegative measurable section of the bundle
$\Omega^{1/2}U=\Omega^{1/2}\ker dt\otimes \Omega^{1/2}\ker ds$.
Using $\rho^U$, we obtain nonnegative sections
$\rho^U_*(f)=(1\otimes \rho^U)(f)$ of the bundle $\Omega^1\ker dt$
and $(\rho^U)^{-1}_*(f)=((\rho^U)^{-1}\otimes 1)(f)$  of the
bundle $\Omega^1\ker ds$; we may then integrate these functions
along the fibers of $t$ and $s$ respectively and obtain functions
$t_!(\rho ^U_{*}(f))$ and $s_!((\rho^U)^{-1}_*f)$ defined on $M$
(with values in $[0,+\infty])$.

\begin{lemma} \label{rhobded}
Thus defined the families $\rho^U$ and $(\rho^U)^{-1}$ are
measurable and bounded over compact subsets of $U$.
\end{lemma}
Note that $\rho^U$ is a section of the one dimensional bundle
$\Hom(\Omega^{1/2}(\ker ds_U), \Omega^{1/2}(\ker dt_U))$, which
explains the above statement.

\begin{proof}
The family $\rho^U$ is continuous on the locally closed subsets of
$U$ on which the dimension of $F_{s(u)}$ is constant; it is
therefore measurable.

We just have to show that $\rho^U$ is bounded in the neighborhood
of every point $u\in U$. Let $u\in U$. Let $\theta $ be a local
diffeomorphism of $M$ carried by $U$. There exists a local
diffeomorphism $\varphi  :U'\to U''$, where $U'$ and $U''$ are
neighborhoods of $u$ in $U$ such that $\varphi(u)=u$, $t_U\circ
\varphi =\theta \circ s_U$ and $s_U\circ \varphi =\theta^{-1}\circ
t_U$. The map $d\varphi $ induces a smooth isomorphism of vector
bundles $\ker ds_U$ and $\ker dt_U$. We just have to compare $\rho
^U$ with the isomorphism $\psi :\Hom(\Omega^{1/2}(\ker ds_U),
\Omega^{1/2}(\ker dt_U))$ induced by $d\varphi $.

Let $v\in U'$; put $x=s_U(v)$ and  $E_v=\ker (ds_U)_v\cap \ker
(dt_U)_v$. The map  $d\varphi _v$ fixes $E_v$ and induces an
automorphism $k_v$ of the vector space $E_v$; the differential
$(d\theta)_x$ induces an isomorphism $\ell_v$ between the
Euclidean vector spaces $F_x$ and $F_{t_U(v)}$. The composition
$(\rho ^U)_v^{-1}\circ \psi_v$ is the multiplication by the square
root of the product $|det (k_v)|.det |\ell_v|$.

Choose a Riemannian metric $g_U$ on $U$, and let $K$ be a compact
subspace of $U'$. Since $d\varphi $ and $d\theta $ are continuous,
they are bounded on $K$ and $s_U(K)$ respectively. It follows that
they multiply the norm of a vector by a ratio bounded by above and
bellow. This gives the desired estimates for $|det (k_v)|.det
|\ell_v|$ when $v\in K$.
\end{proof}

\begin{definition} \label{L1norm}
Let $\cU=(U_i,t_i,s_i)_{i\in I}$ be an atlas.
\begin{enumerate}\renewcommand{\theenumi}{\alph{enumi}}
\renewcommand{\labelenumi}{\theenumi)}
\item Let $(V,t_V,s_V)$ be a bi-submersion adapted to $\cU$. The
$L^1$-norm associated with $\rho^V$ is the map which to $f\in
C_c(V;\Omega^{1/2}V)$ associates $$
\|f\|_{1,\rho^V}=\max\Big(\sup\{ t_{!}(\rho_*|f|)(x),\ x\in
M\},\sup \{s_{!}(\rho_*^{-1}|f|)(x),\ x\in M\}\Big).$$ \item In
the case $(U,t,s)=\coprod_{i\in I} (U_i,t_i,s_i)$, we obtain the
$L^1$-norm $\|f\|_{1,\rho^\cU}=\|f\|_{1,\rho^U}$ for $f\in
C_c(U;\Omega^{1/2}U)=\bigoplus _{i\in I}C_c(U_i;\Omega^{1/2}U_i)$.
\item We note $\|\ \|_{1,\rho}$ the quotient norm in
$\cA_\cU=\bigoplus _{i\in I}C_c(U_i;\Omega^{1/2}U_i)/\cI$.
\end{enumerate}
\end{definition}

It follows from lemma \ref{rhobded} that the quantity
$\|f\|_{1,\rho^V}$ is well defined and finite for every $f\in
C_c(V;\Omega^{1/2}V)$. The map $f\mapsto \|f\|_{1,\rho^V}$ is
obviously a norm, whence  $\|\ \|_{1,\rho}$ is a semi-norm.

\begin{proposition} \label{L1bounded}
The semi-norm $\|\ \|_{1,\rho}$ is a $*$-algebra norm, \ie
satisfies $\|f\star g\|_{1,\rho}\le \|f\|_{1,\rho}\|g\|_{1,\rho}$
and $\|f^*\|_{1,\rho}=\|f\|_{1,\rho}$ (for every $f,g\in
\cA_\cU$).
\end{proposition}
\begin{proof}
Note that in Lemma  \ref{obviouslemma}, one may show that one can
choose $g$ so that $p_!(g)=f$ and $p_!(|g|)=|f|$. It follows that
for every bi-submersion $V$ adapted to $\cU$ the map
$Q_V:C_c(V;\Omega^{1/2}V)\to \cA_\cU$ is norm reducing. Now, it
follows from the Fubini theorem that, if $V,W$ are bi-submersions,
for every $f\in C_c(V;\Omega^{1/2}V)$ and every $g\in
C_c(W;\Omega^{1/2}W)$, we have $\|f\otimes g\|_{1,\rho^{V\circ
W}}\le \|f\|_{1,\rho^V}\|g\|_{1,\rho^W}$. The proposition follows.
\end{proof}

\begin{definition} \begin{enumerate}\renewcommand{\theenumi}{\alph{enumi}}
\renewcommand{\labelenumi}{\theenumi)}
\item We denote by $L_\rho^1(M,\cU,\rho)$ the Hausdorff completion
of $\cA_\cU$ with respect to the $L^1$-norm associated with
$\rho$. When $\cU$ is the path holonomy atlas, we write
$L^1(M,\cF)$ instead of $L_\rho^1(M,\cU)$. \item The {\em full
$C^*$-algebra} of $\cU$ is the enveloping $C^*$-algebra
$C^*(M,\cU)$ of $L_\rho^1(M,\cU)$. Equivalently, it is the
Hausdorff completion of $\cA_\cU$ by the norm $\|(Q_{U})(f)\| =
\sup_{\pi}\{ \|\pi(f)\| \}$ for every $*$-representation $\pi$ of
$\cA_\cU$ on a Hilbert space $\mathcal{H}$ bounded in the sense
$\| \pi(Q_U(f)) \| \leq \| f \|_{1}$. When $\cU$ is the path
holonomy atlas, we write $C^*(M,\cF)$ instead of $C^*(M,\cU)$.
\end{enumerate}\end{definition}

\subsection{Regular representations associated with holonomy pairs and the reduced $C^*$-algebra}

In order to define the reduced $C^*$-algebra of the foliation, we
have to decide which set of representations we consider as being
the regular ones. In fact, there are several possible choices:
\begin{itemize}
\item If the groupoid $G_\cU$ is longitudinally smooth, we can
just consider the corresponding regular representations.
\item In the general case, we may replace $G_\cU$ by a suitable ergodic
decomposition leafwise.
\item We could also decide to consider the representations on square integrable functions on the leaves, or
just consider the (dense set of) leaves without holonomy.
\end{itemize}

In order to give a flexible solution, covering many different
choices, we fix a holonomy pair $(\cU,G)$. In particular, we have a groupoid homomorphism $\alpha:G_\cU\to G$. If $U$ is a bi-submersion adapted to $\cU$, we put $\alpha_U=\alpha\circ q_U:U\to G$.

Fix $x\in M$. The set $G_x=\{\gamma\in G;\ s(\gamma)=x\}$ is a
smooth manifold  by definition of a holonomy pair. We denote by $L^2(G_x)$ the completion of the
space $C_c^\infty(G_x;\Omega^{1/2}_{\C}\ker (ds)_{x})$. This is a
Hilbert space, with inner product $\langle f,g \rangle =
\int_{G_{x}} \Bar{f} \otimes g$ (we can integrate $\Bar{f} \otimes
g$ as a section of $\Omega^{1}(\ker (ds)_{x})$.)

\begin{proposition}
There is a $*$-representation $\pi _x: \cA_\cU \to
\mathcal{L}(L^2(G_x))$ satisfying
$$[\pi_x (Q_{U}(f))(\xi) ](\gamma)= \int_{t_U(u)=t(\gamma)} f(u)\xi(\alpha_U(u)^{-1} \gamma )$$
for every adapted bi-submersion $(U,t_U,s_U)$, every $f \in
C_{c}^{\infty}(U;\Omega^{1/2}U)$ and $\xi \in L^2 (G_x)$. Let
$\rho$ be a positive isomorphism  $\rho :\Omega^{1/2}\ker ds_U\to
\Omega^{1/2}\ker dt_U$ associated  with a Riemannian metric as
above. Then $\pi_x$ extends to a representation of
$L^1(M,\cU;\rho)$.
\end{proposition}

\begin{proof}
To see that $\pi_x$ is well defined, one just needs to check that
if $p:U\to V$ is a morphism of bi-submersions which is a
submersion then $\pi_x (Q_{U}(f))(\xi)= \pi_x
(Q_{V}(p_!(f)))(\xi)$, which is easy.

To show that $\pi_x$ extends to $L^1(M,\cU;\rho)$, we have to show
that for every $f\in C_c^\infty(U;\Omega^{1/2}U)$ and $\xi\in
C_c^\infty(G_x;\Omega^{1/2}G_x)$, we have $\|\pi_x
(Q_{U}(f))(\xi)\|_2\le \|f\|_{1,\rho^U}\|\xi\|_2$.

Endow $U$ and $G$ with the leafwise smooth structures. Then $\alpha_U$
becomes a smooth submersion. For $u,v\in G_x$, put  then
$F(u,v)=(\alpha_{U})_{!}(f)(u^{-1}v)$. Then $F$ is a section of the
bundle of half densities on $G_x\times G_x$. The inequality
follows then from the following classical fact (applied to
$Z=G_x$): Let  $Z$ be a smooth manifold. Fix a positive density
section $\mu$ on $Z$. We thus identify half density sections on
$Z$ and $Z\times Z$ with functions. Let $F$ be a smooth function
on $Z\times Z$. Put  $\|F\|_1=\max\Big(\sup_x \int
|F(x,y)|\,d\mu(y)\,,\,\sup_y \int |F(x,y)|\,d\mu(x)\Big)$. If
$\|F\|_1< \infty$, then the map $\xi\mapsto F\star \xi$ given by
$(F\star \xi )(x)=\int F(x,y)\xi(y)\,d\mu (y)$ is a well defined
bounded map from $L^2(Z)$ into itself with norm $\le \|F\|_1$.

Finally, the fact $\pi_x$ is a $*$-representation follows from the Fubini theorem.
\end{proof}

\begin{definition}
Let $(\cU,G)$ be a holonomy pair. We define the {\em reduced norm}
in $\cA_\cU$ by $\|f \|_{r} = \sup_{x\in M} \|\pi_{x}(f))\|$. The
\emph{reduced $C^*$-algebra} $C_{r}^{*}(G)$ of $G$ is the
Hausdorff-completion of $\cA_\cU$ by the reduced norm. If $\cU$ is
the path holonomy atlas and $G_\cU$ is longitudinally smooth, we
write $C_{r}^{*}(M,\cF)$ instead of $C_{r}^{*}(G_\cU)$.
\end{definition}

\section{Representations of $C^*(M,\cF)$}

In this section, we extend to the (highly singular) groupoid associated
with an atlas $\cU$ of bi-submersions  the following result given
by Renault in \cite{Renault} (see also \cite{FackSkand}): given a locally compact groupoid $G
\gpd M$ there is a correspondence between the representations of
$C^*(G)$ (on a Hilbert space) and the representations of the
groupoid $G$ (on a Hilbert bundle). This  will clarify that the
$L^{1}$-norm we gave earlier provides a good estimate, although
defined up to the choice of a Riemannian metric. As an obvious
consequence, $C^*(M,\cF)$ does not depend on the choice of the
Riemannian metric.

\bigskip
Let us fix again an atlas of bi-submersions $\cU$ for the
foliation $\cF$ on $M$ (not necessarily the path holonomy atlas).
We denote $G=G_{\cU}$ the associated groupoid.  If $U$  is a
bi-submersion adapted to $\cU$, we denote by $q_U:U\to G_\cU$ the
quotient map.

Let $\cA_\cU$ be the associated $*$-algebra. The full
$C^{*}$-algebra will be denoted $C^{*}(\cA_{\cU})$.

\textit{In order to avoid uninteresting technicalities, we consider only representations on separable Hilbert
spaces; we assume that the atlas $\cU=(U_i)_{i\in I}$ is countably
generated, \ie that $I$ is countable and each $U_i$ is
$\sigma$-compact.}

\subsection{Representations of $G_\cU$ and their integration}

We begin with a few considerations about measures.

Let $p : N \to M$ be a submersion and $\mu$ a measure on $M$. Any positive Borel
section $\lambda^p$ of the bundle $\Omega^{1}(\ker dp)$ may be thought
of as a family of positive Borel measures $\lambda^p_{x}$ on $p^{-1}(\{ x
\})$. Such a section defines a measure $\mu \circ \lambda^p$ on $N$ by
$$ \mu\circ\lambda^p(f) = \int_{M}\left(\int_{p^{-1}(\{ x \})}f (y)\,d\lambda^p_{x}(y)\right)d\mu(x) $$
for every $f \in C_{c}(N)$.

\begin{definition}\label{quasi-invariant}
The measure $\mu$ on $M$ is {\em quasi-invariant} if for all
$(U,t_{U},s_{U}) \in \cU$ the measures $\mu\circ\lambda^{t_U}$ and
$\mu\circ\lambda^{s_U}$ are equivalent.
\end{definition}
In other words, there exists a positive locally
$\mu\circ\lambda^{t_U}$-integrable function $D^{U}$ such that
$\mu\circ\lambda^{s_U} = D^{U} \cdot \mu\circ\lambda^{t_U}$.

Let us show that this Radon-Nikodym derivative $D$ gives rise to a
homomorphism $G_\cU\to \R_+^*$ if we make suitable choices of the
densities $\lambda^{t_U}$ and $\lambda^{s_U}$.

We fix a Riemannian metric on $M$. We saw in 4.4 that it induces a
Borel isomorphism $\rho^{U} :\Omega^{1}(\ker ds_{U}) \to
\Omega^{1}(\ker dt_{U})$ for every bi-submersion $(U,t_{U},s_{U})
\in \mathcal{U}$. Choose a section $\lambda^{s_U} \in
C_{c}(U;\Omega^{1}(\ker ds_{U}))$ and take the corresponding
$\lambda^{t_U} = \rho^{U}(\lambda^{s_U})$ which is a positive
Borel section, bounded over compact sets of the bundle
$\Omega^{1}(\ker dt_{U})$. Let $\mu \circ \lambda^{s_U}$ and
$\mu\circ\lambda^{t_U}$ be the associated measures of $U$.
\begin{enumerate}
\item The function $D^{U}$ depends on the chosen Riemannian metric but does not depend on the choice of
  $\lambda^{s_U}$. This follows from the linearity of the map
  $\rho^{U}$ applied to the product of $\lambda^{s_U}$ with a
positive continuous function.
\item The maps $D^{U}$ may be replaced by a map $D : G_{\cU} \to \R^{+}$
  defined by $D(q_{U}(u)) = D^{U}(u)$ for any bi-submersion of the
  atlas $\cU$. More precisely,

Let $\varphi :U\to V$ be a morphism of bi-submersions. Then $D^{U}
= D^{V} \circ \phi$ almost everywhere.

Indeed, the statement is local: we may replace $U$ by a small open
neighborhood of a point $u$. Thanks to corollary
\ref{corsmoothlong}, we may assume that $\varphi$ is a submersion.
Finally, we may assume that $U=L\times V$ where $L$ is a manifold
and $\varphi$ is the projection. We may then take $\lambda
^{s_U}=\lambda^L\times \lambda^{s_V}$ where $\lambda^L $ is a $1$
density on $L$; then $\lambda ^{t_U}=\lambda^L\times
\lambda^{t_V}$, and the statement follows.

\item To show that $D$ is a {\em homomorphism} we need to examine
the behavior of the $D^{U}$s with respect to the composition of
bi-submersions. For two bi-submersions $(U_{1},t_{1},s_{1})$ and
$(U_{2},t_{2},s_{2})$ in $\cU$ denote $(U_{1} \circ U_{2},t,s)$
their composition. Then $s^{-1}(x) = U_{1}
\times_{s_{1},t_{1}}(U_{2})_{x}$ for every $x \in M$. Notice that
$\ker (ds)_{(u_{1},u_{2})}$ contains $\ker (ds_1)_{u_{1}} \times \{
0 \}$ and the quotient is isomorphic to $\ker (ds_2)_{u_{2}}$. So
we obtain  a canonical isomorphism
$$\Omega^{1}(\ker (ds)_{(u_{1},u_{2})} = \Omega^{1}(\ker (ds_1)_{u_{1}})
\otimes \Omega^{1}(\ker (ds_2)_{u_{2}})$$
Likewise, we find
$$\Omega^{1}(\ker (dt)_{(u_{1},u_{2})} = \Omega^{1}(\ker (dt_1)_{u_{1}})
\otimes \Omega^{1}(\ker (dt_2)_{u_{2}})$$
One checks that, up to these isomorphisms, we have an equality $\rho_{U_{1}\circ U_{2}} =
\rho_{U_{1}} \otimes \rho_{U_{2}}$. An easy computation gives then
$$D^{U_{1}\circ U_{2}} = D^{U_{1}} \circ p_{1} \cdot D^{U_{2}} \circ p_{2},$$
where $p_{1}, p_{2}$ are the projections of $U_{1} \circ U_{2}$ to
$U_{1}$ and $U_{2}$ respectively.
\end{enumerate}

Due to the last two observations we may think of $D$ as being a homomorphism $D :
  G_{\cU} \to \R^{+}$.

In the same way, we define representations of $G_\cU$ as being
`measurable homomorphisms from $G_\cU$ to the unitary group.' The
precise definition is the following:

\begin{definition} \label{defrep}
A {\em representation} of $G_{\mathcal{U}}$ is a triple $(\mu, H, \pi)$ where:
\begin{enumerate}\renewcommand{\theenumi}{\alph{enumi}}
\renewcommand{\labelenumi}{\theenumi)}
\item $\mu$ is a quasi-invariant measure on $M$; \item $H =
(H_{x})_{x \in M}$ is a measurable (with respect to $\mu$) field
of Hilbert spaces over $M$. \item For every bi-submersion
$(U,t,s)$ adapted to $\cU$,  $\pi^U$ is a measurable (with respect
to $\mu\circ \lambda $) section of the field of unitaries
$\pi^U_u:H_{s(u)}\to H_{s(u)}$;
\end{enumerate}

Moreover, we assume:
\begin{enumerate}
\item $\pi $ is `defined on $G_\cU$':\\
if $f:U\to V$ is a morphism of bi-submersions, for almost all
$u\in U$ we have $\pi_u^U=\pi_{f(u)}^V$.
\item $\pi $ is a homomorphism: \\
If $U$ and $V$ are bi-submersions adapted to $\cU$, we have
$\pi^{U\circ V}_{(u,v)}=\pi^U_{u}\pi^V_{v}$ for almost all
$(u,v)\in U\circ V$.
\end{enumerate}
\end{definition}

Any representation $(\mu, H, \pi)$ where $\mu$ is quasi-invariant
gives rise to a representation of $C^{*}(M,\cF)$ on the space
$\cH=\int_M^\oplus H_x\, d\mu(x)$  of sections of the Hilbert
bundle $H$:

For every bi-submersion $U$ adapted to $\cU$, define
$\hat{\pi}_U : C_{c}(U,\Omega^{1/2}U) \to B(\cH)$  by putting
$$\hat{\pi}_U(f)(\xi)(x) = \int_{U^{x}}(1\otimes\rho^{U}) (f(u))\pi^U_{u}(\xi(s_{U}(u))) D^U(u)^{1/2}$$
$\mu$-a.e. for every $f \in C_{c}(U,\Omega^{1/2}(U)),\  \xi \in
\cH$ and $x \in M$. 

One checks that we thus define a representation of $C^*(M;\cF)$ in a few steps \begin{enumerate}
\item If $\varphi : U\to V$ is a morphism  of bi-submersions which is a submersion, for every $f \in C_{c}(U,\Omega^{1/2}(U))$, we obviously have $\hat{\pi}_U(f)=\hat{\pi}_V(\varphi_!(f))$.

\item Let $(U_{1},t_{1},s_{1})$ and
$(U_{2},t_{2},s_{2})$ be two bi-submersions in $\cU$, and denote $(U_{1} \circ U_{2},t,s)$ their composition. For every $f_1\in C_{c}(U_1,\Omega^{1/2}(U_1))$ and $f_2 \in C_{c}(U_2,\Omega^{1/2}(U_2))$, we have $\hat{\pi}_{U_1}(f_1)\hat{\pi}_{U_2}(f_2)=\hat{\pi}_{U_1\circ U_2}(f_1\otimes f_2)$.

\item Let $(U,t,s)$ be a bisubmersion in $\cU$, and let $U^{-1}=(U,s,t)$ be its inverse. Denote by $\kappa:U\to U^{-1}$ the identity. We have $\hat \pi_{U^{-1}}(\overline f\circ \kappa)=\big(\hat{\pi}_U(f)\big)^*$.

It follows from 1. that there is a linear map $\hat \pi:\cA_\cU\to \cL(\cH)$ such that for every $U$ adapted to $\cU$ we have $\hat\pi_U=\hat\pi\circ Q_U$; it then follows from 2. and 3. that the map $\hat \pi $ is a $*$-homomorphism.

\item The proof of \cite{Renault} can be adapted to show that $\hat{\pi}$ is continuous with
respect to the $L^{1}$-norm (\cf definition \ref{L1norm}), and therefore defines a representation $\hat \pi$ of $C^*(M,\cU)$ on $\cH$.
\end{enumerate}

\subsection {Desintegration of representations}
We now show that every representation of $C^*(M,\cU)$ is of the
above form \ie it is the integrated form of a representation of
the groupoid $G_\cU$. Actually, the proof given by J. Renault in
\cite{Renault} can be easily adapted to our case. We give here an
alternate route, based on Hilbert $C^*$-modules.

Let us first note that $C_0(M)$ sits in the multipliers of
$C^*(M,\cU)$: if $(U,t_U,s_U)$ is a bi-submersion adapted to
$\cU$, for $f\in C_c(U;\Omega^{1/2})$ and $h\in C_0(M)$, we just
put $h.Q_U(f)=Q_U((h\circ t)f)$ and $Q_U(f).h=Q_U(f(h\circ s))$.
It follows that every non degenerate representation of
$C^*(M,\cU)$ on a separable Hilbert space $\cH$ gives rise to a
representation of $C_0(M)$ and therefore a measure (class) $\mu $
on $M$ and a measurable field of Hilbert spaces $(H_x)_{x\in M}$
such that $\cH=\displaystyle\int^\oplus \!H_x\,d\mu(x)$ on which
$C_0(M)$ is naturally represented. We are going to show that $\mu$
is quasi-invariant and $G_\cU$ is represented on the Hilbert field
$(H_x)_{x\in M}$.

Recall (\cf \eg \cite{Connes-Skandalis}) that along with a
submersion  $p:N\to M$ is naturally associated a Hilbert $C_0(N) -
C_0(M)$-bimodule $\cE_p$  which is the continuous family of
Hilbert spaces $\big(L^2(p^{-1}(\{x\})\big)_{x\in M}$, on which
$C_0(N)$ is represented by multiplication. More specifically,
$\cE_p$ is the completion of $C_c(N;\Omega^{1/2}\ker dp)$ with
respect to the $C_0(M)$ valued scalar product given by the formula
$\langle \xi,\eta\rangle(x)=\int_{p^{-1}(\{x\})} \overline \xi
\eta$. The operation of $C_0(N)$ by multiplication obviously
extends to $\cE_p$.

Now, if $p:N\to M$ is a submersion, we put $\tilde \cE
_p=\cE_p\otimes _{C_0(M)}C^*(M,\cU)$. A typical element of $\tilde
\cE_p$ is the class of a section $C_c(N;\Omega^{1/2}\ker
dp)\otimes C_c(U;\Omega^{1/2})$ where $U$ is a bi-submersion
adapted to $U$; we therefore associate elements of $\tilde \cE_p$
to every $\varphi\in C_c(N\times U;\Omega^{1/2}(\ker dp\times \ker
dt_U\times \ker ds_U))$; it is easily seen that the image of such
a $\varphi $ only depends on its restriction to $N\times
_{p,t_U}U$. Now, $N\times _{p,t_U}U$ is easily seen to be a
bi-submersion between $(M,\cF)$ and $(N,p^{-1}\cF)$. It is in fact
quite easy to prove:

\begin{proposition} \label{prop5.4}
For every bi-submersion $(V,t_V,s_V)$ between $(M,\cF)$ and
$(N,p^{-1}\cF)$ such that $(V,p\circ t_V,s_V)$ is adapted to
$\cU$, we have a map $Q_V:C_c(V,\Omega^{1/2}V)\to \tilde \cE_p$.
The images of $Q_V$ span a dense subspace of $\tilde \cE_p$.  If
$V,W$ are two bi-submersions as above between $(M,\cF)$ and
$(N,p^{-1}\cF)$, $f\in C_c(V,\Omega^{1/2}V)$, $g\in
C_c(W,\Omega^{1/2}W)$ and $h\in C_0(N)$, we have the formulae
$$hQ_V(f)=Q_V((h\circ t_V)f)\ \ \ \hbox{and}\ \ \ \langle
Q_V(f),Q_{W}(g)\rangle = Q_{V^{-1}\circ W}(\overline f\otimes
g).\eqno(1)$$
\end{proposition}

Let us give a few explanations on the statement of this proposition:
\begin{itemize}
\item $\Omega^{1/2}V$ denotes the bundle $\Omega^{1/2}(\ker
dt_V\times \ker ds_V)$; \item $V^{-1}\circ W=V\times
_{t_V,t_{W}}W$; it is a bi-submersion of $(M,\cF)$; \item we
identify $\cA_\cU$ with its image in $C^*(M,\cU)$, thus
$Q_{V^{-1}\circ W}(\overline f\otimes g)\in C^*(M,\cU)$.
\end{itemize}

\begin{proof}
One checks easily formulae (1) for bi-submersions of the form
$N\times _{p,t_U}U$ and $f,g$ of the form $f_1\otimes f_2$ where
$U$ is a bi-submersion of $(M,\cF)$, $f_1\in C_c(N,\Omega^{1/2})$,
$f_2\in C_c(U,\Omega^{1/2})$. In general, $V$ is covered by open
subsets of this form, and the conclusion follows using partitions
of the identity.
\end{proof}

\begin{remnot} \label{remnott}
Let $p,q:N\to M$ be two submersions such that
$p^{-1}\cF=q^{-1}\cF$; assume moreover that for every
bi-submersion $(V,t_V,s_V)$ of $(M,\cF)$ between $(M,\cF)$ and
$(N,p^{-1}\cF)$, $(V,p\circ t_V,s_V)$ is adapted to $\cU$ if and
only if $(V,q\circ t_V,s_V)$ is adapted to $\cU$. It follows from
proposition \ref{prop5.4} that there is a canonical isomorphism
between $\sigma _p^q\in \cL(\tilde \cE_p,\tilde \cE_q)$.
\begin{enumerate}
\renewcommand{\theenumi}{\alph{enumi}}
\renewcommand{\labelenumi}{\theenumi)}
\item Obviously, $\sigma_p^p=\id$, $\sigma_q^p=(\sigma_p^q)^{-1}$;
furthermore, if $r:N\to M$ is a third submersion satisfying the
same properties, we have $\sigma _p^r=\sigma _q^r\sigma _p^q$.

\item Let $(U,t_U,s_U)$ be a bi-submersion adapted with $\cU$. To
simplify notation, we will just denote by $\sigma^U\in \cL(\tilde
\cE_{s_U},\tilde \cE_{t_U})$ the unitary operator
$\sigma_{s_U}^{t_U}$.
\end{enumerate}
\end{remnot}

\medskip
Now let $\varpi :C^*(M,\cU)\to \cL(\cH)$ be a representation of
$C^*(M,\cU)$. Since $C_0(M)$ sits in the multiplier algebra of
$C^*(M,\cU)$, the representation $\varpi $ gives rise to a
representation of $C_0(M)$. We thus get a measure (class) $\mu $
on $M$ together with a measurable family of Hilbert spaces
$(H_x)_{x\in M}$. If $p:N\to M$ is a submersion, write
$\cE_{p}\otimes _{C_0(M)}\cH=\int ^{\oplus}_N H_{p(y)}\, d\mu\circ
\lambda(y)$: this gives the representations of $C_0(N)$ on
$\cE_{p}\otimes _{C_0(M)}\cH$. Let $(U,t_U,s_U)$ be a
bi-submersion adapted to $\cU$. The representations of $C_0(U)$ on
$\cE_{s_U}\otimes _{C_0(M)}\cH=\tilde \cE_{s_U}\otimes
_{C^*(M,\cU)}\cH$ and $\cE_{t_U}\otimes _{C_0(M)}\cH=\tilde
\cE_{t_U}\otimes _{C^*(M,\cU)}\cH$ are isomorphic through the
unitary operator $\alpha^U=\sigma^U\otimes
_{C^*(M,\cU)}\id_{\cH}$. It follows that $\mu$ is quasi-invariant;
furthermore the operator $\sigma^U\otimes _{C^*(M,\cU)}\id_{\cH}$
yields a measurable family $(\pi^U_u)_{u\in U}$ where
$\pi^U_u:H_{s_U}\to H_{t_U}$ is a unitary operator.

One sees easily that:
\begin{enumerate}
\item If $f:U\to V$ is a morphism of bi-submersions,
$\pi^U_u=\pi^V_{f(u)}$ (almost everywhere):
\begin{enumerate}
\item this is obvious if $U$ is an open subset of $V$, \item this
is also easy if $f:U=L\times V\to V$ is the projection, \item from
these facts the formula is established if $f$ is a submersion;
\item the general case follows from corollary \ref{corsmoothlong}.
\end{enumerate}
\item Let $(U,t_U,s_U)$ and $(V,t_V,s_V)$ be bi-submersions. Let $(W,t_W,s_W)$ be the
bi-submersion $U\circ V$. We have $\pi_{(u,v)}^W=\pi_u^U\pi_v^V$ (almost everywhere in $W$).\\
Indeed, let $\alpha:W\to M$ be the map $(u,v)\mapsto
s_U(u)=t_V(v)$. The map $(u,v)\mapsto u$ is a morphism of
bi-submersions between $(W,t_W,\alpha)$ and $U$ and the map
$(u,v)\mapsto v$ is a morphism of bi-submersions between
$(W,\alpha, s_W)$ and $V$. The equality follows from
\ref{remnott}.
\end{enumerate}

It follows from the above that $(\mu,H,\pi)$ is a representation
of the groupoid $G_\cU$ in the sense of definition \ref{defrep}.

Finally, one checks that $\varpi$ is the representation $\hat \pi$
associated with $\pi$.

\section{Further developments}

We end with a couple of remarks that we will expand elsewhere:

\subsection{The tangent groupoid, longitudinal pseudo-differential operators and analytic index}

\subsubsection {The cotangent space}

The $K$-theory class of a symbol should take place in the
$K$-theory of the total space of a cotangent bundle. The cotangent
bundle of the foliation is the ``bundle'' $(\cF^*_x)_{x\in M}$.
Let us discuss this space.

\begin{definition}
Every $\xi \in\cF_x^*$ with $x\in M$ defines a linear functional
on $\cF$ as a composition $\xi\circ e_x:\cF\to \cF_x\to \R$. Let
$\cF^*$ be the union of $\cF_x^*$ endowed with the topology of
pointwise convergence on $\cF$. It is a locally compact space.
\end{definition}

\subsubsection{The tangent groupoid}
Let $\lambda :M\times \R\to\R$ be the second projection. Consider
the foliation $\cF\times \{0\}$ of $M\times \R$. Denote by $T\cF$
the foliation $\{\lambda X;\ X\in \cF\times \{0\}\}$.

The groupoid of this foliation is $\bigcup_{x\in M}\cF_x\times
\{0\}\cup G\times \R^*$. It is called the \emph{tangent groupoid.}
Its $C^*$-algebra contains as an ideal $C_0(\R^*)\otimes
C^*(M;\cF)$ with quotient $C_0(\cF^*)$.

This tangent groupoid allows to construct a map $K(C_0(\cF^*))\to
K(C^*(M;\cF))$. This map is the \emph{analytic index}.

\subsubsection{Pseudo-differential calculus along the foliation}

One can also define the longitudinal pseudo-differential operators
associated with a foliation. These are unbounded multipliers of
the $C^*$-algebra.

The differential operators are very easily defined: They are
generated by vector fields in $\cF$.

Let $(U,t,s)$ be a bi-submersion and $V\subset U$ an identity
bisection. Denote by $N$ the normal bundle to $V$ in $U$ and let
$a$ be a (classical) symbol on $N^*$. Let $\chi$ be a smooth
function on $U$ supported on a tubular neighborhood of $V$ in $U$
and let $\phi:U\to N$ be an inverse of the exponential map (defined
on the neighborhood of $V$). A pseudodifferential kernel on $U$ is
a (generalized) function $k_a:x\mapsto \int a(x,\xi)
\exp(i\phi(u)\xi) \chi (u)\,d\xi$ (the integral is an oscillatory
integral, taken over the vector space $N^*_x$).

As in the case of foliations and Lie groupoids (\cf \cite{Connes0,
Monthubert-Pierrot, Nistor-Weinstein-Xu}), we show:

\begin{itemize} \item The kernel $k_a$ defines a multiplier of $\cA(\cF)$.
\item Those multipliers when $U$ runs over an atlas, $V$ covers
$M$ and $a$ runs over the classical symbols form an algebra. This
algebra only depends on the class of the atlas.

\item The algebra of pseudodifferential operators is filtered by
the order of $a$. The class of $k_a$ only depends up to lower
order on the restriction of the principal part of $a$ on $\cF$.

\item Negative order pseudodifferential operators are elements of
the $C^*$-algebra (both full and reduced) of the foliation.  Zero
order pseudodifferential operators define bounded multipliers of
the $C^*$-algebra of the foliation.

We have an exact sequence of $C^*$-algebras $$0\to C^*(M,\cF)\to
\Psi^*(M,\cF)\to C_0(S^*\cF)\to 0$$ where $\Psi^*(M,\cF)$ denotes
the closure of the algebra of zero order pseudodifferential
operators.

\item  One can also take coefficients on a smooth vector bundle over $M$.

\item Elliptic operators of positive order (\ie operators whose
principal symbol is invertible when restricted to $\cF$) give rise
to regular quasi-invertible operators  (\cf \cite{Vassout}).
\end{itemize}

\subsubsection{The analytic index}

Elliptic (pseudo)-differential operators have an index which is an
element of $K_0(C^*(M,\cF))$. One easily sees that this index
equals the one given by the `tangent groupoid' and the
corresponding exact sequence of $C^*$-algebras.

\subsection{Generalization to a `continuous family' case}

As it was explained by Paterson in \cite{Paterson}, one doesn't
really need a Lie groupoid in order to have a nice
pseudodifferential calculus: the only thing which matters is the
fact that the $s$ and $t$  fibers are smooth manifolds. {\tt Examples of longitudinaly smooth groupoids naturally appear also in the case of stratified Lie groupoids recently studied in a fundamental paper by Fernandez, Ortega and Ratiu (\cite{FOR}).}

In the same way, we wish to consider singular foliations on
locally compact spaces, which are only smooth in the leaf
direction. Actually, one way to define those foliations is to
start with an atlas of bi-submersions. Indeed, in the smooth case,
an atlas of bi-submersions associated with $\cF$, obviously
defines the foliation $\cF$.

We may then define the holonomy groupoid, the $C^*$-algebra
exactly as above, together with a suitable pseudodifferential
calculus.

\end{document}